\documentclass[hidelinks,11pt]{amsart}
\usepackage{amsmath,amsthm,amscd,amssymb,amsfonts}
\usepackage{geometry}
\usepackage{subcaption}
\usepackage{titlesec}

\usepackage{pb-diagram}
\usepackage{graphicx}
\usepackage{setspace}
\usepackage{array}
\usepackage{pstricks}
\usepackage{pst-node}
\usepackage{url}
\usepackage{xcolor}
\usepackage{soul}
\usepackage{hyperref}
\usepackage{ulem}
\usepackage[most]{tcolorbox}
\usepackage{float}
\usepackage{multicol}
\usepackage{mathrsfs}
\usepackage{MnSymbol}
\usepackage{gensymb} 
\usepackage{enumerate}
\usepackage{mathtools} 
\usepackage[abs]{overpic}

\newtheorem{thm}{Theorem}[section]

\theoremstyle{definition}

\newtheorem{defn}[thm]{Definition}

\newtheorem{ques}[thm]{Question}

\numberwithin{equation}{section}
\numberwithin{figure}{section}

\usepackage{titlesec}
\titleformat{\section} {\normalfont\bfseries\filcenter}{\thesection }{1em}{}
\titleformat{\subsection} {\normalfont\bfseries\filcenter}{\thesubsection}{1em}{}
\renewcommand \thesubsection{}

\begin{document}

\title[{Topological Data Analysis of Mortality Patterns During the COVID-19 Pandemic}]{Topological Data Analysis of Mortality Patterns During the COVID-19 Pandemic}

\author[MEGAN FAIRCHILD]{MEGAN FAIRCHILD}
\address{The University of Utah}
\email{megan.fairchild@utah.edu}

\author[MATTHEW LEMOINE]{MATTHEW LEMOINE}
\address{Department of Mathematics \\ Louisiana State University}
\email{mlemo36@lsu.edu}

\begin{abstract}
    Topological Data Analysis is a relatively new field of study that uses topological invariants to study the shape of data. We analyze a dataset provided by the Centers for Disease Control and Prevention (CDC) using persistent homology and MAPPER. This dataset tracks mortality week-to-week from January 2020 to September 2023 in the United States during the COVID-19 pandemic. We examine the dataset as a whole and break the United States into geographic regions to analyze the overall shape of the data. Then, to explain this shape, we discuss events around the time of the pandemic and how they contribute to the observed patterns.
\end{abstract}

\maketitle

\section{Introduction} \label{sec:intro}

Topological data analysis (TDA) is a powerful tool for analyzing the geometric structure of a dataset. The dataset must be finite, and we must have a notion of \textit{distance} between two points. We then \textit{build a continuous structure} using the points of the dataset that has some underlying topology. This structure is called a \textit{simplicial complex} and can be thought of as a graph (vertices and edges) with higher dimensional objects (faces, tetrahedrons, etc.). We then extract information from the simplicial complex using the main tool called \textit{persistent homology}. This allows us to identify interesting features of our dataset and analyze the structure of the data from a topological standpoint.  Useful information about TDA and persistent homology were first explained in Carlsson's work in \cite{Carlsson}.

A given dataset may live in higher dimensions than we can visualize, so having a way to project the data points down to a lower dimension gives us a way to analyze clusters and get an idea of what to look at in our persistent homology. One such dimension reduction technique is called MAPPER. This technique allows us to use  principal component analysis to look at our dataset as a graph in 2 or 3 dimensions. A reader interested in learning about the MAPPER algorithm in glorious detail is referred to Madukpe, Ugoala, and Zulkepli \cite{reviewmapperalgorithm}.

In this work, we analyze data focused on the mortality rates in the US before, during, and after the COVID-19 pandemic. We use topological data analysis to infer information about the shape of our data. There have been several works discussing topological data analysis as it relates to the pandemic, we refer the interested reader to \cite{covid_mapper}, \cite{phang_nonlinear_2024}, and \cite{so_topological_2021}. In these works, they use topological data analysis to infer information about the spread of COVID-19. Specifically in \cite{phang_nonlinear_2024}, they use a dataset of daily infections from the beginning of the pandemic to June 2024 to look at the spread of COVID-19 in Malaysia. Recently, Assaf, Rammal, Goupil, Kacim, and Vrabie discussed how to reduce the number of false positives from COVID-19 using topological data analysis on images of the lungs of a patient \cite{assaf_topological_2025}. In contrast to other work that has been done, we look more into the aftermath of the pandemic by looking at mortality rates in various geographical regions in the United States. One goal of this work is to determine how the spread of the coronavirus affected mortality rates in the US and to see if topological data analysis can serve as an effective tool to answer this question.

The dataset we focus on in this paper consists of weekly counts of deaths and causes of death from January 2020 to September 2023, as reported by the Centers for Disease Control and Prevention (hereafter referred to as ``the CDC”) on their open source data platform, see \cite{CDC}. The data is categorized by region into five groups: West, South, Northeast, Midwest, and outlying territories. Using these regional distinctions, we applied two tools from topological data analysis to gain insight into the overall structure and shape of the dataset. We aim to answer the following questions. 

\begin{ques}\label{question1}
    What are the main contributing factors to the overall shape of the dataset for the United States versus the geographical regions? 
\end{ques}

\begin{ques}\label{question2}
    Are the topological features depicted in the barcode analysis also depicted in the MAPPER projection?
\end{ques}

\begin{ques}\label{question3}
    How large of an impact did the COVID-19 mortality rate have on the overall shape of the various datasets?
\end{ques}

We begin by reviewing necessary background in topology and topological data analysis in Section \ref{sec:background}. In Section \ref{sec:methods and data}, we describe our dataset obtained from the CDC and outline our methods. Section \ref{sec:barcode} presents the analysis of the barcodes generated from the data. In Section \ref{sec:Mapper}, we discuss the application of the MAPPER algorithm and the associated visualization tools, and we use the MAPPER outputs to test our hypotheses. Also in the MAPPER section, we summarize the results from our data interpretation from both the barcodes and the MAPPER, and Section \ref{sec:conclusions} presents the overall conclusions and final observations regarding the shape and structure of the dataset. We include all MAPPER figures and projections in Appendix \ref{sec:Mapper_figs}. In the MAPPER analysis section \ref{sec:Mapper}, we only include one of the figures to highlight specific features of a given region.

\textbf{Acknowledgments.} The authors were partially supported by NSF Research Training Groups in the Mathematical Sciences (RTG) Grant No. NSF-DMS 2231492. 

\section{Background} \label{sec:background}

We wish to compute persistent homology of our dataset, so we must discuss how we can construct the shape of a dataset. Following Carlsson and Silva \cite{CarlssonSilva} and Adams and Tausz \cite{javaplex}, we define an $n$-simplex as the \textit{convex hull} of $n$ vertices, which is the intersection of all convex sets containing the vertices. An abstract \textit{simplicial complex}, $X$, is defined by a set $Z$ of vertices, or $0$-simplices, and for each $k\geq 1$, we have a set of $k$-simplices $\sigma = [z_0,z_1,...,z_k]$, where $z_i \in Z$. Each $k$-simplex has $k+1$ faces obtained by deleting one of the vertices. The following inclusion property must be satisfied: if $\sigma \in X$ then for all $\delta \subset \sigma$, $\delta \in X$. Geometrically, this is ensuring that if a $k$-simplex is in our complex that all $(k-1)$-simplices in the $k$-simplex must also be in the complex. Additionally, the intersection of any two simplices within the simplicial complex must be \textit{closed.} That is, the intersection of two simplices is again a simplex.

In the geometry of our objects, we think of $0$-simplices as vertices, $1$-simplices as edges, $2$-simplices as triangular faces, and $3$-simplices as solid tetrahedrons. We then study the Betti numbers $b_{k}$, which in a more geometric sense, measure the number of $k$-dimensional holes in a topological space. In particular, $b_0$ measures the number of connected components. 

A \textit{filtration} of a simplicial complex $X$ is a collection of subcomplexes defined with a parameter $t$, $\{X(t)\ |\ t \in \mathbb{R}\}$, where each $X(t)$ is a subcomplex of $X$ and $X(t) \subseteq X(t')$ whenever $t\leq t'$. The \textit{filtration value} of a simplex $\sigma$ is the smallest $t$ such that $\sigma \in X(t)$. We will often use the term \textit{stream} to refer to a filtered simplicial complex.

Having a filtration (or stream) gives us inclusion maps from each step in our filtration to the next. These inclusion maps induce a chain maps in homology. Then we can take homology to get the persistent homology using these maps by looking for holes or voids that persist from one step in the filtration to another, see \cite{Hatcher}, \cite{virk}, and \cite{WeiWei} for more information about homology and persistent homology.

\subsection{Vietoris-Rips Streams}

In order to compute persistent homology of a dataset in any meaningful way, we must have a notion of distance between two points in the dataset. Let $d$ denote the distance, or metric, on our point cloud dataset $Z$. 

\begin{defn}[Vietoris-Rips Stream \cite{javaplex}]

    The Vietoris-Rips stream has a complex $VR(Z, t)$ for some $t$ defined as follows. 
    \begin{enumerate}
        \item The vertex set is $Z$.
        \item given $x, y \in Z$, the edge $(x, y)$ is included in $VR(Z, t)$ if $d(x, y) \leq t$.
        \item a higher dimensional simplex is included in $VR(Z, t) $ if all of its faces are. 
    \end{enumerate}
    
\end{defn}

Note that this defines a filtered complex with filtration given by $t$. To build some intuition, consider a point cloud in Euclidean space equipped with the standard metric. Around each point, we draw an $\epsilon$-ball for $\epsilon>0$. We let $\epsilon$ grow continuously, and record it's change in distance with the perspective of time $t$. The reader can imagine that as the $\epsilon$-balls intersect, simplices are added to the complex: two intersecting balls define an edge, three define a face, and so on. When we begin discussing analysis of our dataset, we will use the terms `seconds' or `time $t$' to refer to the $t$ value of our filtration.

\section{Methods and Data} \label{sec:methods and data}

For this project, we used the GUDHI Python library to construct the simplicial complex, compute persistent homology, and generate barcodes. We also applied the MAPPER algorithm, as further described in Section \ref{sec:Mapper}, to view both 2D and 3D projections of our dataset. The dataset was obtained from the CDC via one of their open-source platforms \cite{CDC} and spans January 2020 through September 2023, recording weekly mortality counts for each state. The dataset columns include week, year, state, and cause of death. A detailed analysis of the barcodes is presented in Section \ref{sec:barcode}. In Section \ref{sec:conclusions}, we compare the barcode results with the MAPPER outputs to provide a cohesive analysis of the dataset.

In addition to analyzing data for the entire United States (hereafter referred to as “Whole US”), we examined mortality patterns by geographical region. The country was divided into five regions, as outlined in Table \ref{tab:Regions}, to allow for more detailed regional analysis. Each region contains at least 1,000 entries, with the exception of the Non-Contiguous US, which includes approximately 580 entries. As a side note, we have included Washington D.C., Maryland, and Delaware in the Northeast to keep the regions as close to the same size as possible. According to the U.S. Census Bureau \cite{census}, these three states and city are included in the South.

\begin{table}[ht]
    \centering
    \begin{tabular}{|l|l|}
        \hline
        Region/Group & States in this Region (number of states)\\ \hline 
        West & California, Arizona, New Mexico, Nevada, Utah, Colorado, Wyoming, \\ & Oregon, Idaho, Washington, Montana (11)\\[10pt] \hline
        Midwest & Missouri, Kansas, Illinois, Indiana, Ohio, Nebraska, Iowa, Michigan, \\ & South Dakota, Wisconsin, Minnesota, North Dakota (12) \\[10pt]  \hline
        Northeast & Washington D.C. , Maryland, Delaware, New Jersey, Pennsylvania, \\ & New York City, New York,  Connecticut, Rhode Island, \\ &  Massachusetts, Vermont, New Hampshire, Maine (13) \\[10pt] \hline
        South & Texas, Florida, Louisiana, Mississippi, Alabama, Georgia, \\ & South Carolina, Arkansas, Oklahoma, Tennessee, North Carolina, \\ & Kentucky, Virginia, West Virginia (14)\\[10pt] \hline
        Non-Contiguous US & Alaska, Hawaii, Puerto Rico (3)\\[10pt]
        \hline
    \end{tabular}
    \caption{Regions and Groups of States}
    \label{tab:Regions}
\end{table}

\subsection{Code details and parameters}
Our code is available on GitHub \cite{our_code} for those interested. To begin our analysis, each column of the dataset was normalized. We then constructed the Vietoris-Rips complex using the standard Euclidean metric and a filtration parameter \texttt{max\_edge\_length = 2.0}. Note that we chose \texttt{max\_edge\_length = 2.0} as increasing the parameter adds more simplices to the complex and it became too computationally expensive. We then decided our maximum dimension of a simplex allowed in our complex is \texttt{max\_dimension = 2}, again this was determined based on the issue with expensive computation due to the size of the dataset. The higher dimensional analysis of the barcodes was done, if we thought that there could be a chance of a higher dimensional feature. This is discussed as it comes up in the barcode analysis in section \ref{sec:barcode}. Note that \texttt{max\_dimension = 2} means we will look at dimensions 0 and 1 in our barcode analyses.

Next, we construct a simplex tree, which is an efficient data structure from the GUDHI library used to represent the filtration. Each node in the simplex tree represents a filtration, and it stores the filtration value, which is the smallest scale at which the simplex appears. This allows us to traverse the tree quickly and query simplices and their filtration values. We additionally track which simplices contribute to each persistent feature and export this information to a CSV file. Essentially, we track which components are merging when a new edge appears in the filtration. 

Finally, we plot the persistence barcodes. In the barcode diagram, each bar represents a topological feature and spans the interval of its existence in the filtration. The diagrams allow us to view the births and deaths of relevant topological features and view the persistent homology of the simplicial complex on our dataset. The barcode analysis is discussed in great detail in Section \ref{sec:barcode}.

One of the most prominent tools in TDA is the Mapper algorithm, developed by Singh, M\'emoli and Carlsson \cite{Mapper}, which provides a graphical summary of the data's topological structure. We use the \texttt{KeplerMapper} Python library and explore the impact of parameter choices in both dimensionality reduction and clustering on the resulting topological summaries. To facilitate the construction of the Mapper graph and reduce computational complexity, we apply Principal Component Analysis (PCA), a linear dimensionality reduction technique. PCA helps remove noise while preserving the most important structure in the data. It does this by identifying new axes (called principal components) along which the data varies the most. These components are ranked by the amount of variance they capture — that is, how much the data spreads out along each direction — allowing us to retain only the most informative features. We choose to reduce the dataset to 2 and 3 dimensions, respectfully, and analyze both scenarios. The PCA-reduced then serves as the input for clustering and Mapper construction.

Density-Based Spatial Clustering of Applications with Noise (DBSCAN) is then applied to the reduced data. The choice of clustering algorithm is crucial for Mapper since it defines the nodes in the resulting simplicial complex. DBSCAN requires two parameters, epsilon, \texttt{eps}, and minimum number of samples, \texttt{min\_samples}. Epsilon is distance between points, and the minimum number of samples is the minimum number of points required to form a dense region (i.e., a cluster). A point is considered ``core" if it has at least \texttt{min\_samples} points (including itself) within its \texttt{eps}-radius. We configure DBSCAN with $\texttt{eps}=30$ and $\texttt{min\_samples}=10$. Finally, we implement the projection into 2-dimensional space, or 3-dimensional space, using the MAPPER algorithm. Further details on the MAPPER algorithms and images are provided in Section \ref{sec:Mapper}.

We took great care in analyzing appropriate parameters for accuracy in our analysis to ensure the integrity of the geometry of the dataset between MAPPER and Vietoris-Rips constructions. We compared pairwise distances in the original high-dimensional space and the PCA-reduced space. Specifically, we computed the full pairwise distance matrices using the Euclidean distance. This comparison is crucial because both the Mapper algorithm (via DBSCAN clustering) and Vietoris–Rips complexes rely heavily on pairwise distances to identify neighborhoods, clusters, or simplices. If dimensionality reduction significantly distorts these distances, the topological features extracted may not reflect the true shape of the original data. This analysis was motivated by concerns that topological constructions like the Vietoris–Rips complex and clustering-based Mapper graphs might produce inconsistent representations of the data’s shape. By comparing the distance ranges, we gain confidence that the major geometric relationships in the data are preserved during projection.

\section{Barcode Analysis}\label{sec:barcode}

In this section, we will examine the barcode for each individual region, followed by an analysis of the barcode for the entire United States. We will explore possible structures within our dataset. In these barcodes, each bar represents a feature’s lifespan—indicating when it emerged (was born) and when it disappeared (died). The red bars in our figures are telling us about the persistence in the $0$-th homology group (i.e. the connected components). The blue bars are telling us about the persistence in the first homology group (i.e. $S^1$ features that persist). For example, if you wanted to know how many connected components there were in our filtration at time $t = 1$, you would draw a line slicing the bars at $t = 1$ and count the number of red bars.

As previously mentioned in Section \ref{sec:methods and data}, our maximum dimension of a simplex in our simplicial complex is 2. Thus, we analyze only the $0^{th}$ and $1^{st}$ homology groups. Our maximum filtration value was also set to 2, which is referenced as the variable $t$ for time in this section. Some regions we set the maximum $t$ value to 1.75 to better view relevant homological features.

\subsection{Overarching Notes}\label{subsection:overarching_notes}

The first notable observation is that the barcodes vary significantly across regions. These regions were grouped based on their geographic location, yet the shape of each dataset differs from one region to another. Additionally, a notable feature in the data that included week and year is the dramatic increase in the number of $S^{1}$ components around time $t\approx 1$. This corresponds to the time scale of the data, where both year and week progress is in discrete increments of 1.

In some cases, the barcodes that include the weeks and years looks similar, if not the same, as the barcodes that exclude the weeks and years (e.g. The South). This tells us that the weeks and years did not play a significant role in the shape of our data. 

Moving forward for brevity, we will say a dataset that includes weeks and years is a dataset ``with dates'', and a dataset that does not include weeks nor years is a dataset ``without dates''. We are interested in comparing and contrasting these datasets for each region to determine how much weight the dates columns are contributing to the overall topology of our dataset. For example, the South dataset, analyzed below, is a great example of a region where the dates columns do not appear to play a significant role.

\subsection{Whole US}\label{sub:wholeusbarcodes}

The Whole US is the data taken by summing up each entry in every state for every week. So the Whole US has the fewest entries at about 200. Below are the barcodes for the Whole US including and excluding dates.

\begin{figure}[H]
	\begin{subfigure}[b]{0.49\textwidth}
		\includegraphics[width=\textwidth]{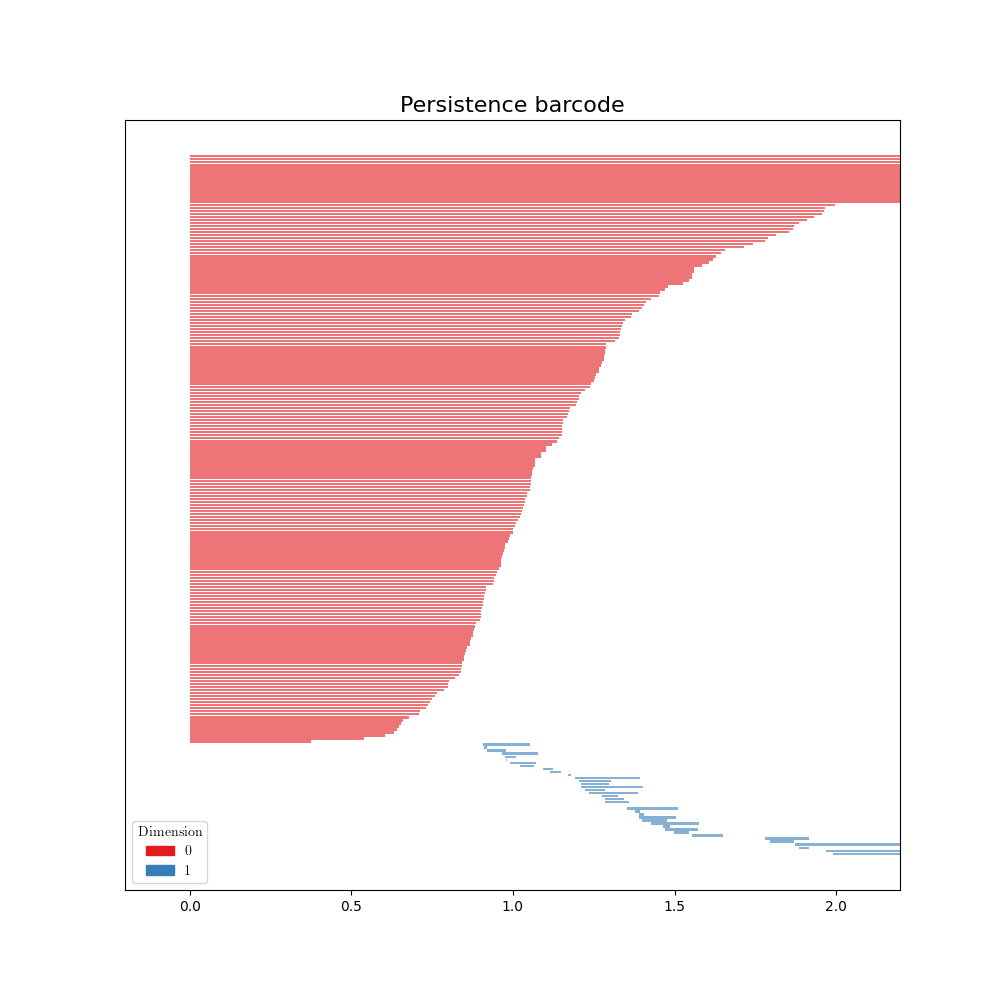}\vspace{-11pt}
                \caption{The Whole US (with dates)}\label{fig:US_bars_dates}
	\end{subfigure}
	~
	\begin{subfigure}[b]{0.49\textwidth}
		\includegraphics[width=\textwidth]{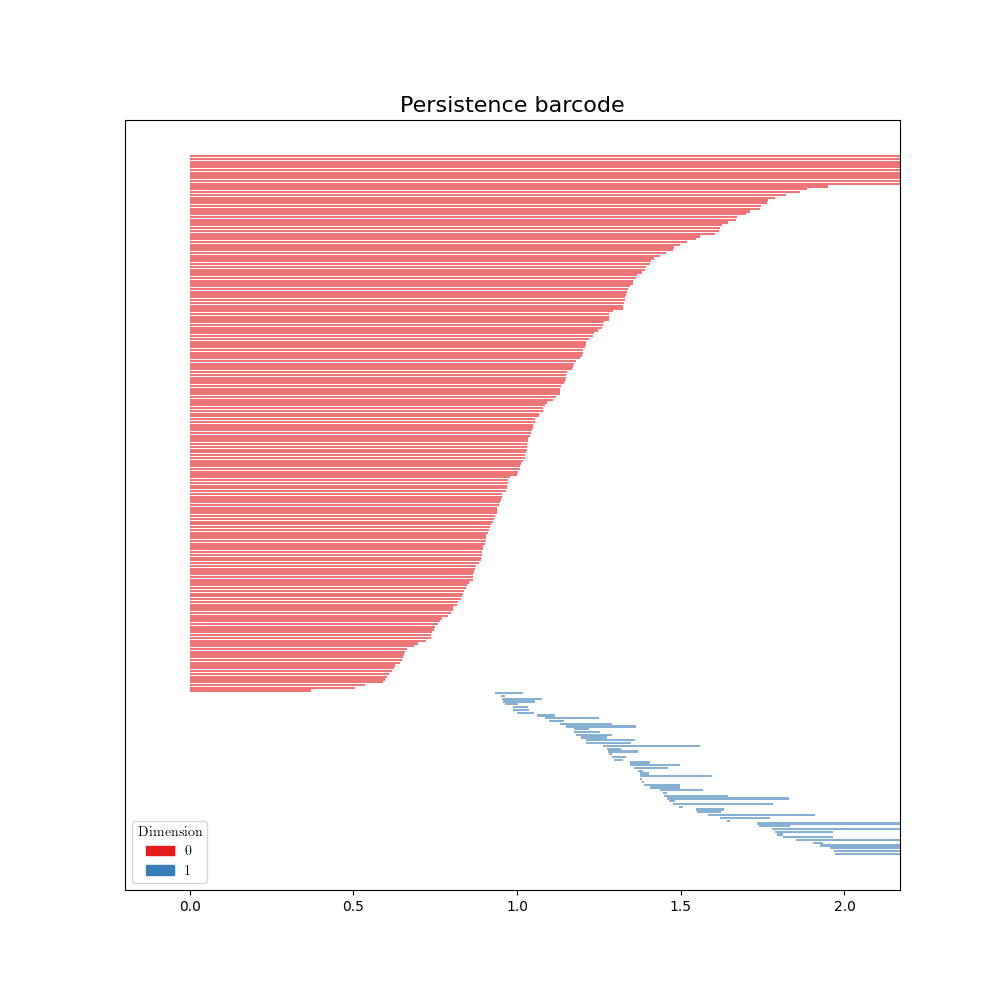}\vspace{-11pt}
                \caption{The Whole US (without dates)}\label{fig:US_bars_nodates}
	\end{subfigure}
    \caption{The Whole US barcodes, where the $x$-axis is our time.}
    \label{fig:US_bars}
\end{figure}

We begin by noting that our barcodes in the figures above look very similar. This tells us that whether we are looking at the Whole US with the dates or without the dates we get similar information. We again see a smooth joining pattern like we saw with the south dataset. Again, this is telling us that the data points are distributed in a normal way. We see less $S^1$ components in both of these figures. This is due to the limited number of entries that we have here. And we do not see the dramatic increase in the $S^1$ components that we have seen in all the regional datasets. This is because when we look at the country as a whole the difference of weeks and years gets smoothed out and we do not see the dramatic consequence of this on the bigger scale. In this dataset, we conclude that our data looks like a point cloud with a few arms sticking out causing the $S^1$'s to show up, but not a major topological shape.

\subsection{West}\label{sub:westbarcodes}

The West, as was defined in Section \ref{sec:methods and data}, consists of the following states: California, Arizona, New Mexico, Nevada, Utah, Colorado, Wyoming,
Oregon, Idaho, Washington, and Montana, for a total of 11 states. Below are the barcodes for the West with and without dates.

\begin{figure}[H]
	\begin{subfigure}[b]{0.49\textwidth}
		\includegraphics[width=\textwidth]{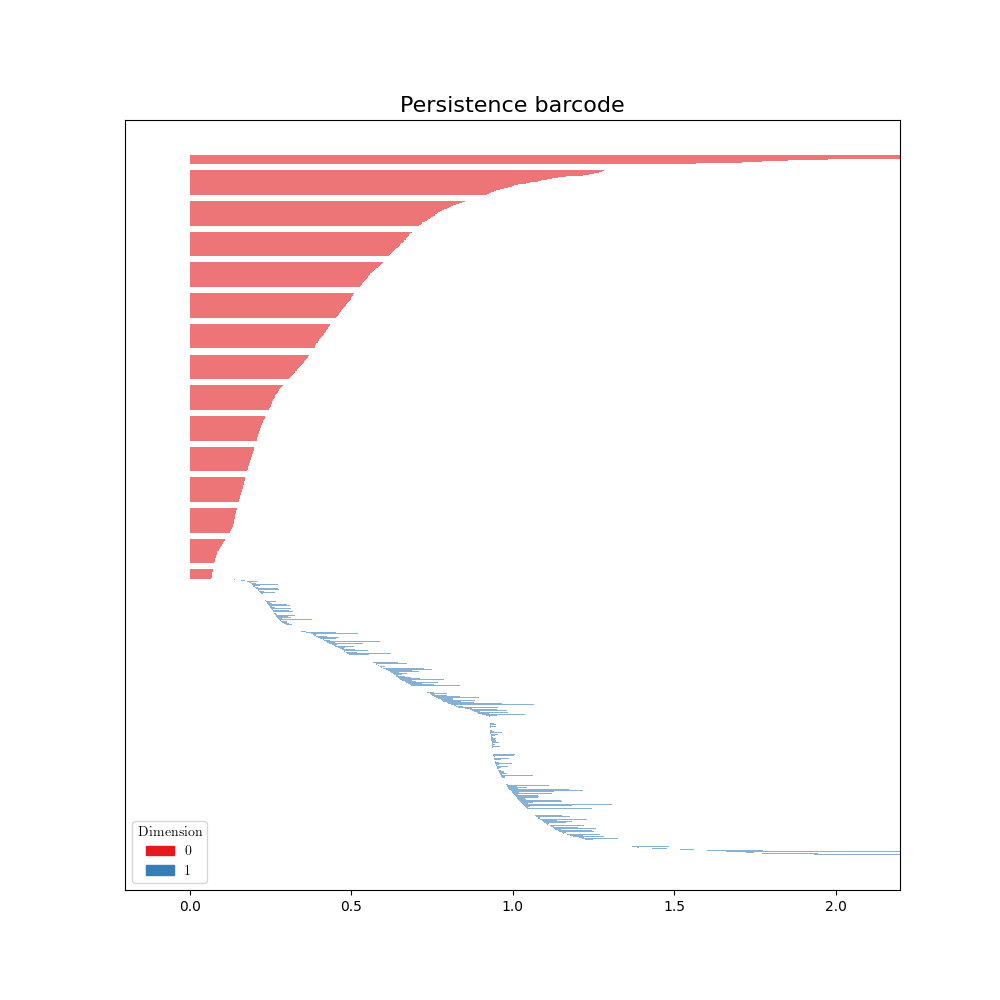}\vspace{-11pt}
                \caption{West (with dates)}\label{fig:west_bars_dates}
	\end{subfigure}
	~
	\begin{subfigure}[b]{0.49\textwidth}
		\includegraphics[width=\textwidth]{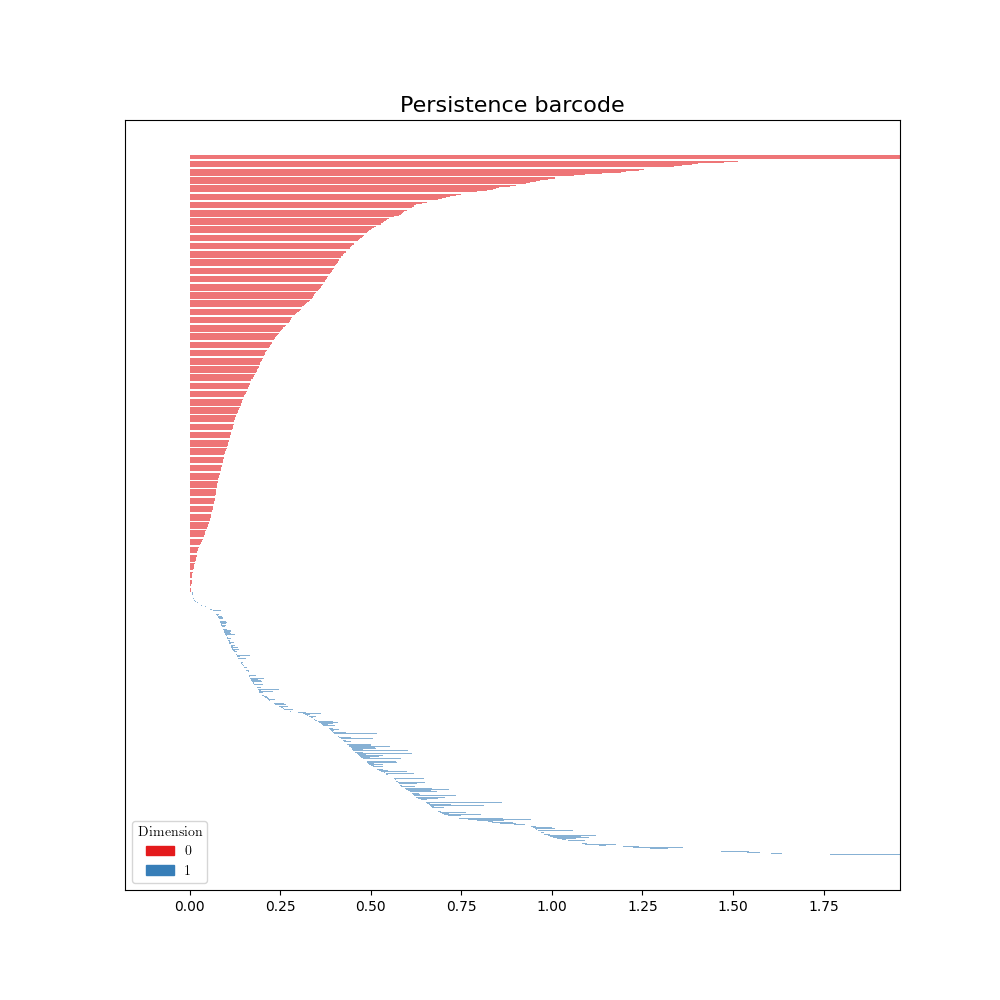}\vspace{-11pt}
                \caption{West (without dates)}\label{fig:west_bars_nodates}
	\end{subfigure}
    \caption{The West barcodes, where the $x$-axis is our time.}
    \label{fig:west_bars}
\end{figure}

Now that we have pointed out the prominent features in this barcode, we can begin to figure out what these features tell us about the shape of the data. 

We begin our discussion of prominent features in Figure \ref{fig:west_bars_dates}. We see there is a dramatic increase in $b_{1}$ at $t \approx 1$. This feature in our barcode comes from the weeks and years, giving extra integer values to each of our entries. This is artificially adding an extra partition to our dataset that is separated by distance 1 for each year and each week. We also notice that in both Figures \ref{fig:west_bars_dates} and \ref{fig:west_bars_nodates}, the disconnected components merge quickly into one connected component. By $t \approx 1.25$, we have one connected component. This feature is not unique to this barcode, and we will see it in other barcodes, but it is still an important feature to note. This tells us that all the data points are `close' to each other, with few outliers. We begin to answer Question \ref{question3} by noticing that COVID-19 mortality rates did not cause a completely separated dataset. 

Next, in Figure \ref{fig:west_bars_nodates}, we see that the 1-dimensional holes are short-lived. That is to say, the $S^{1}$'s do not persist for very long. This leads us to conclude that the points in this dataset are tightly clustered with only a few outliers. We also observe that although the data points are relatively close to one another, they form small, tightly packed clusters or ``circles" that are small in scale compared to the entire structure.

When imagining what this structure may look like, we imagine a radio tower with small clusters around it. Upon running the MAPPER algorithm, we see a very similar shape (Figure \ref{fig:west_3D_noDates}), further cementing our hypothesis is true, and partially answering Question \ref{question2}.

\subsection{Midwest}\label{sub:midwestbarcodes}

The Midwest, as was defined in Table \ref{tab:Regions}, comprises the following states: Missouri, Kansas, Illinois, Indiana, Ohio, Nebraska, Iowa, Michigan, South Dakota, Wisconsin, Minnesota, and North Dakota, for a total of 12 states. Below are the barcodes that we obtained from the Midwest data with and without the dates included, respectfully. 

\begin{figure}[H]
	\begin{subfigure}[b]{0.49\textwidth}
		\includegraphics[width=\textwidth]{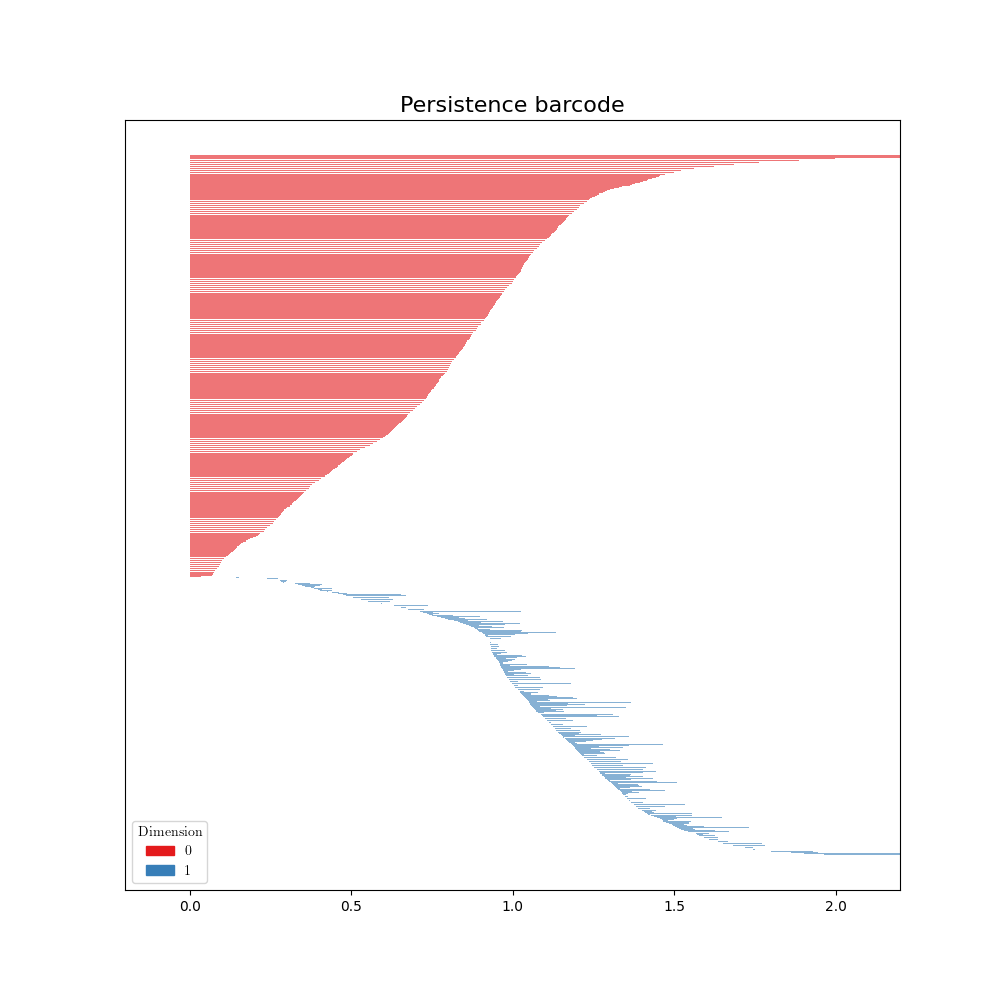}\vspace{-11pt}
                \caption{Midwest (with dates)}\label{fig:midwest_bars_dates}
	\end{subfigure}
	~
	\begin{subfigure}[b]{0.49\textwidth}
		\includegraphics[width=\textwidth]{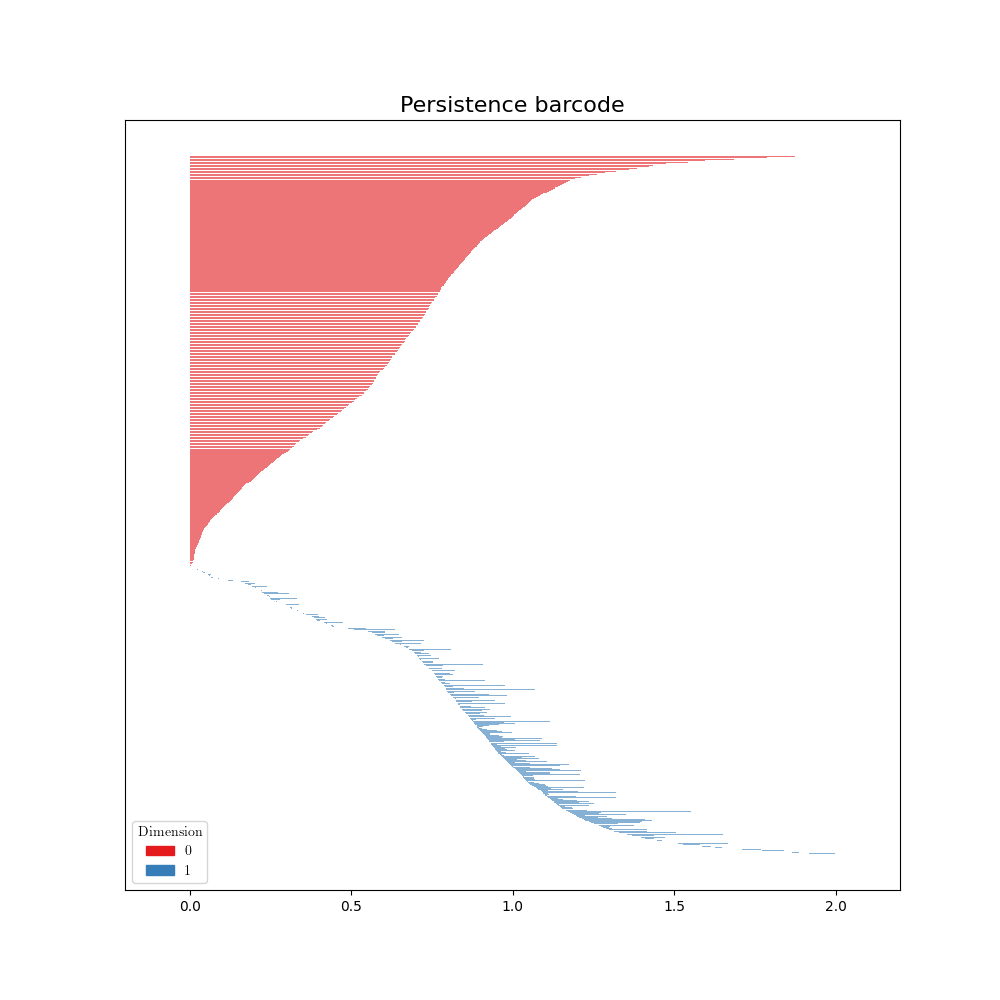}\vspace{-11pt}
                \caption{Midwest (without dates)}\label{fig:midwest_bars_nodates}
	\end{subfigure}
    \caption{The Midwest barcodes, where the $x$-axis is our time.}
    \label{fig:midwest_bars}
\end{figure}

Looking at these barcodes in figures \ref{fig:midwest_bars_dates} and \ref{fig:midwest_bars_nodates}, we see they look similar, which tells us that the including the dates in our dataset does not contribute to the shape of our dataset. We again see the increase in the $S^1$ components in figure \ref{fig:midwest_bars_dates} that we saw before. And again this comes from the years and weeks being included, and they all differ by 1. 

Now we note that there is a steady joining of the $b_0$ components until $t\approx 1.4$. This is telling us that the data points are all spaced out in a linear way. Meaning that as our epsilon balls are increasing the rate at which our components are joining up is the same. So our points are all spaced out following a linear pattern. This means that either we have several clusters that are all spaced out this way or there is one cluster that has this spacing pattern. There being multiple clusters would explain the bump in the $b_0$ components at $t\approx 0.75$ and the resulting $b_1$ components that come after this. We again see that most of the components are merged together by $t\approx 1.5$. The fact that most of the components are joined together by then means that there are very few outliers in this dataset. 

We conclude that our dataset for the Midwest, does not have a shape that depends on the dates. We also conclude that this dataset looks like a point cloud where the points are spaced out linearly from a `center' or there are a few clusters that are spaced out in a linear way.

\subsection{Northeast}\label{sub:northeastbarcodes}

The Northeast, as was defined in section \ref{sec:methods and data}, is comprised of the following states and cities: Washington D.C. , Maryland, Delaware, New Jersey, Pennsylvania, New York City, New York, Connecticut, Rhode Island, Massachusetts, Vermont, New Hampshire, and Maine, for a total of 13 states or cities. Below are the barcodes that we obtained from the code for the Northeast with and without the dates included.

\begin{figure}[H]
	\begin{subfigure}[b]{0.49\textwidth}
		\includegraphics[width=\textwidth]{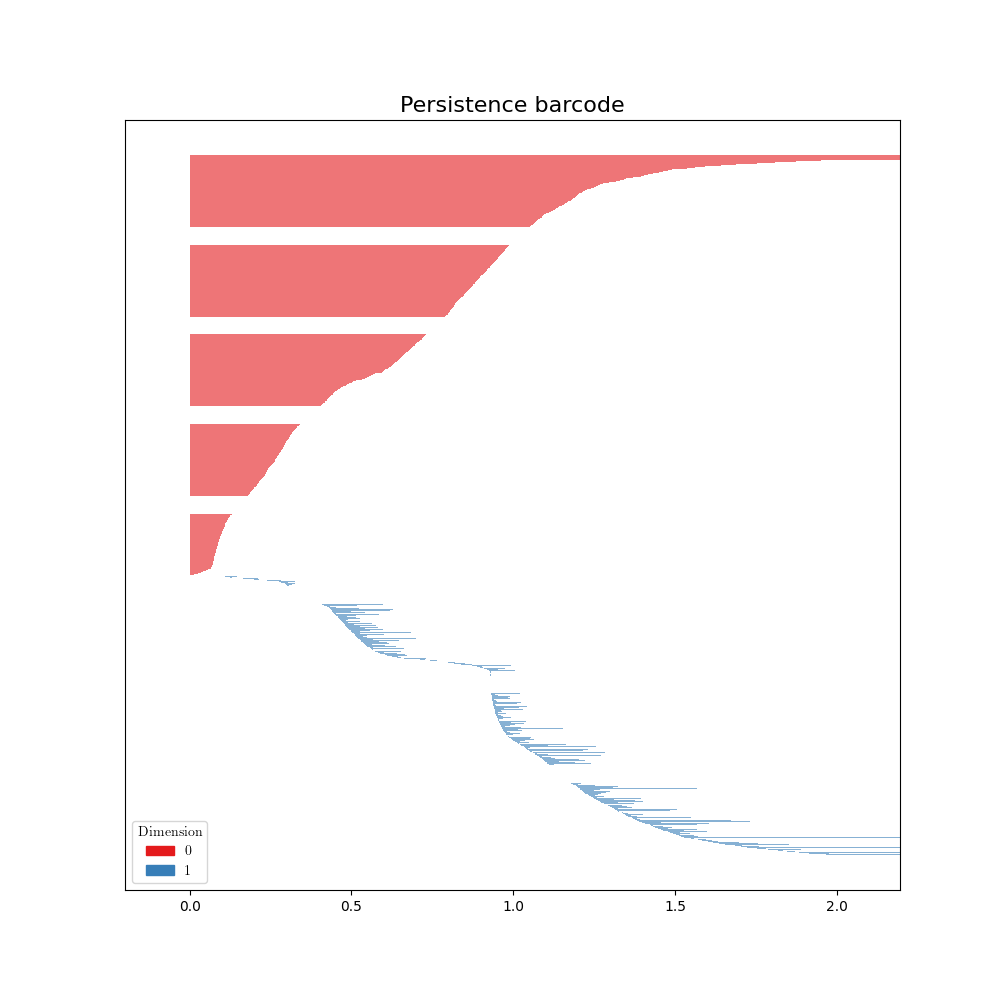}\vspace{-11pt}
                \caption{Northeast (with dates)}\label{fig:northeast_bars_dates}

	\end{subfigure}
	~
	\begin{subfigure}[b]{0.49\textwidth}
		\includegraphics[width=\textwidth]{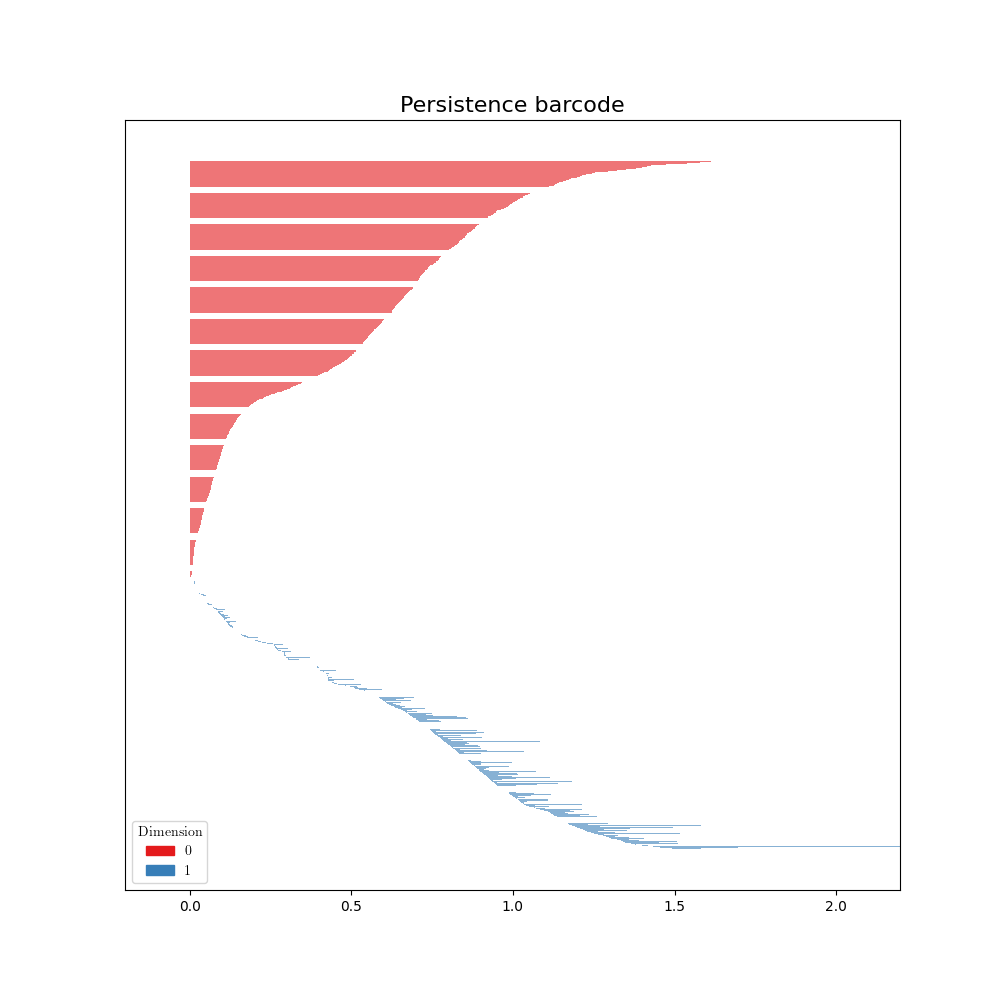}\vspace{-11pt}
                \caption{Northeast (without dates)}\label{fig:northeast_bars_nodates}
	\end{subfigure}
    \caption{The Northeast barcodes, where the $x$-axis is our time.}
    \label{fig:northeast_bars}
\end{figure}

Looking at the key features of these barcodes, we see first that there is a dramatic bump in the $b_0$ components at $t\approx 0.5$. This feature is less pronounced in the barcode including that dates than excluding the dates. We see this bump because there are clusters forming early in the stream ($t <0.4$) that get merged together at $t\approx 0.5$. These clusters are less pronounced when the dates are included, which means that the dates break up these clusters. Also in the $b_0$ components we see that most of our components are merged together by $t\approx 1.75$. Therefore, our clusters take care of most of the outlying points and there aren't any outlying clusters.

Again, we see the increase in the $S^1$ components at $t\approx 1$ in figure \ref{fig:northeast_bars_dates} which is coming from the weeks and years. In the $S^1$ components in figure \ref{fig:northeast_bars_dates}, we also see that these last for longer than in figure \ref{fig:northeast_bars_nodates}. This is also coming from the weeks and years. We don't see these long bars in figure \ref{fig:northeast_bars_nodates}.

Therefore, we can make some conclusions about the shape of our datasets. The dataset that includes dates closely resembles a point cloud that has some striation to it. There is a separation in the clusters that is caused by the weeks and years being included. In the dataset that does not include the dates, we see that these clusters that were separated in the previous dataset are actually apart of the same cluster. These could have come from a spike in cases during a particular time in the year that repeats every year. Thus for the dataset that excludes the dates, we have a point cloud that is comprised of a few clusters.

\subsection{South}

The South, as was defined in section \ref{sec:methods and data}, comprises the following states: Texas, Florida, Louisiana, Mississippi, Alabama, Georgia, South Carolina, Arkansas, Oklahoma, Tennessee, North Carolina, Kentucky, Virginia, and West Virginia, for a total of 14 states. Below are the barcodes we obtained for the South when we included and excluded dates.

\begin{figure}[H]
	\begin{subfigure}[b]{0.49\textwidth}
		\includegraphics[width=\textwidth]{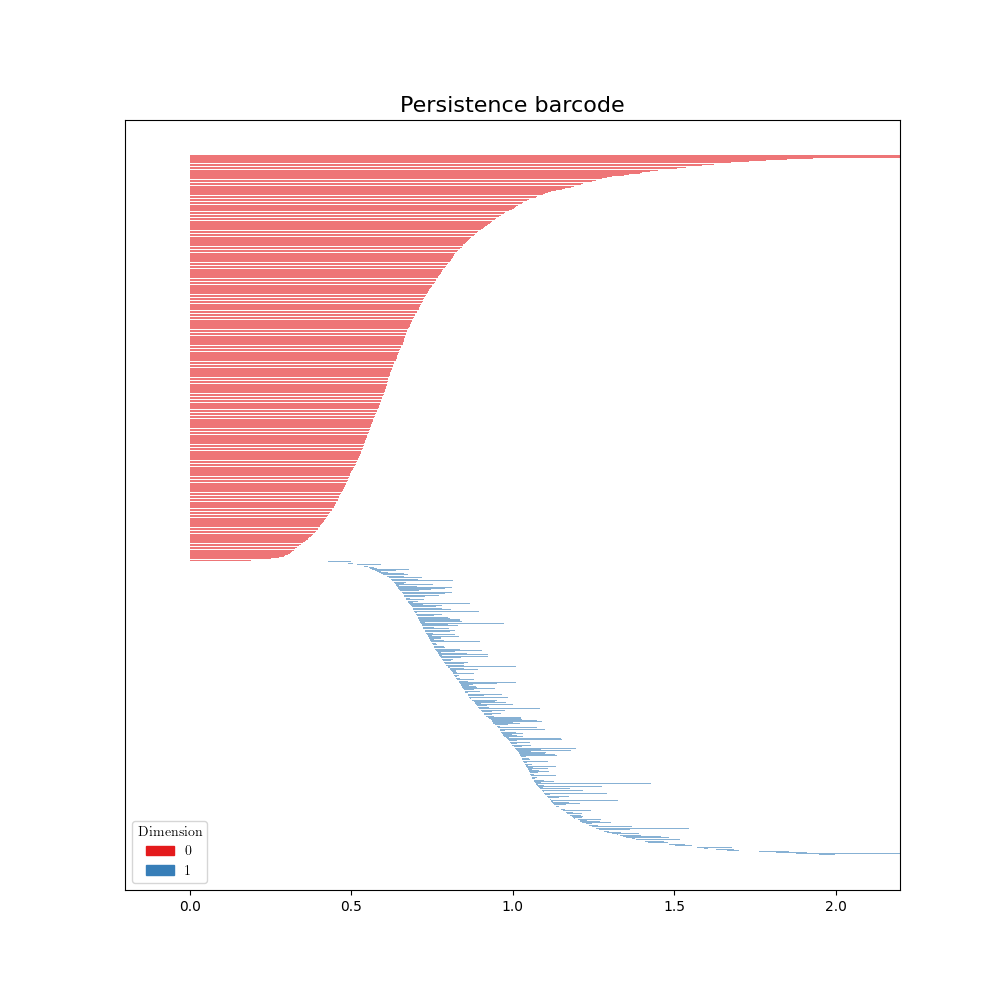}\vspace{-11pt}
                \caption{South (with dates)}\label{fig:south_bars_dates}
	\end{subfigure}
	~
	\begin{subfigure}[b]{0.49\textwidth}
		\includegraphics[width=\textwidth]{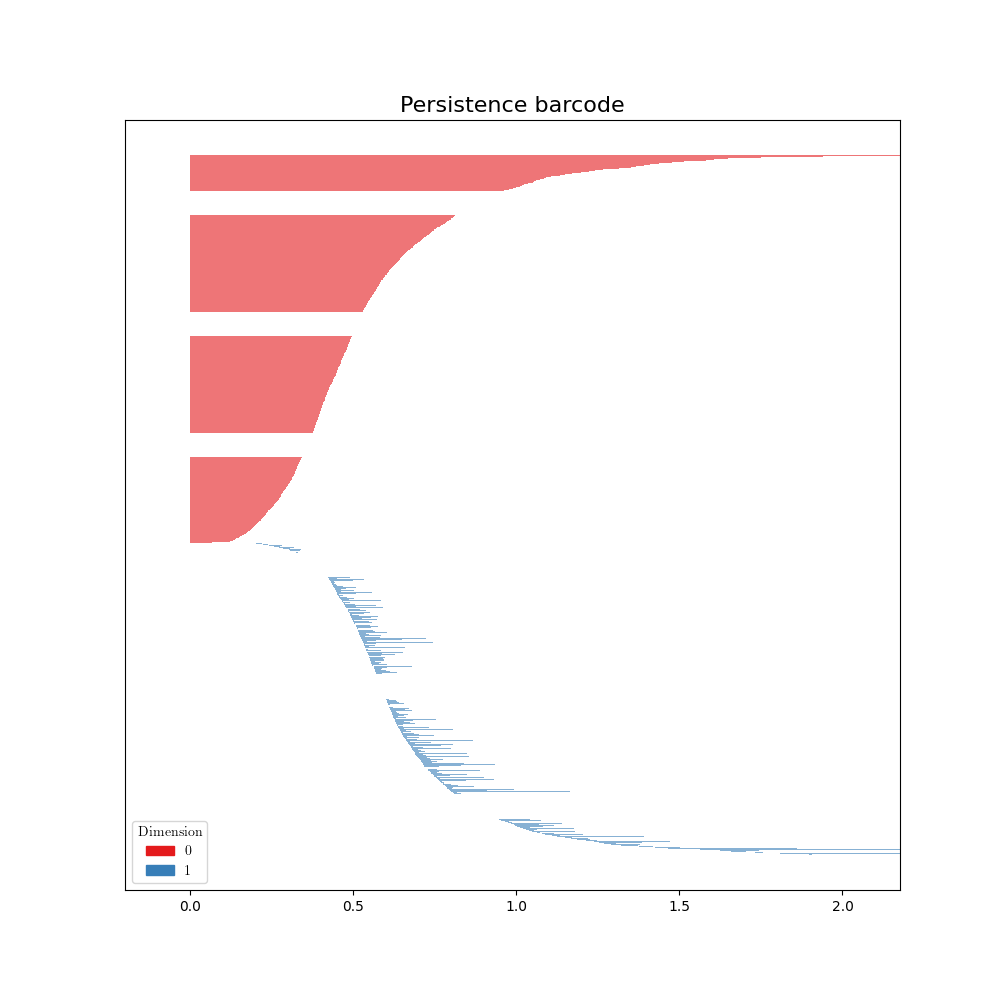}\vspace{-11pt}
                \caption{South (without dates)}\label{fig:south_bars_nodates}
	\end{subfigure}
    \caption{The South barcodes, where the $x$-axis is our time.}\label{fig:South_bars}
\end{figure}

First we note that these bars look very different from our other barcodes. We also note that these barcodes look very similar. We see smooth merging patterns for the $b_0$ components. This pattern of merging is telling us that our points are spaced out in a normal way. Meaning that our points are normally distributed. We can see this because as our points merge together, they merge quicker and quicker, then they taper off and merge less quickly. The pattern of distributed points also explains our $S^1$'s not lasting for long amounts of time. The smoothness of our $b_0$ components joining together would seem to imply that we do not have any auxiliary clusters, except one major cluster. But we see some trailing off in both of these barcodes after $t\approx 1$. This is telling us that there are some outliers that are still normally distributed away but are outliers nonetheless. 

Also we note that these barcodes are very similar, with one notable difference being there is a slight bump in the $S^1$ components at $t\approx 1$ in figure \ref{fig:south_bars_dates}. This is again coming from the weeks and years being in increments of 1. But this bump is not very large, which tells us that the dates do not play a huge role in the shape of our dataset. We conclude that our dataset looks like a large point cloud where the points are distributed normally over our space with some `center', and these points make few if any clusters.

\subsection{Non-Contiguous US}\label{sub:noncontiguousbarcodes}

The Non-Contiguous US, as was defined in section \ref{sec:methods and data}, comprises the following states and territories: Alaska, Hawaii, and Puerto Rico, for a total of 3 states and territories. Below are the barcodes obtained from these areas when we included and excluded dates.

\begin{figure}[H]
	\begin{subfigure}[b]{0.49\textwidth}
		\includegraphics[width=\textwidth]{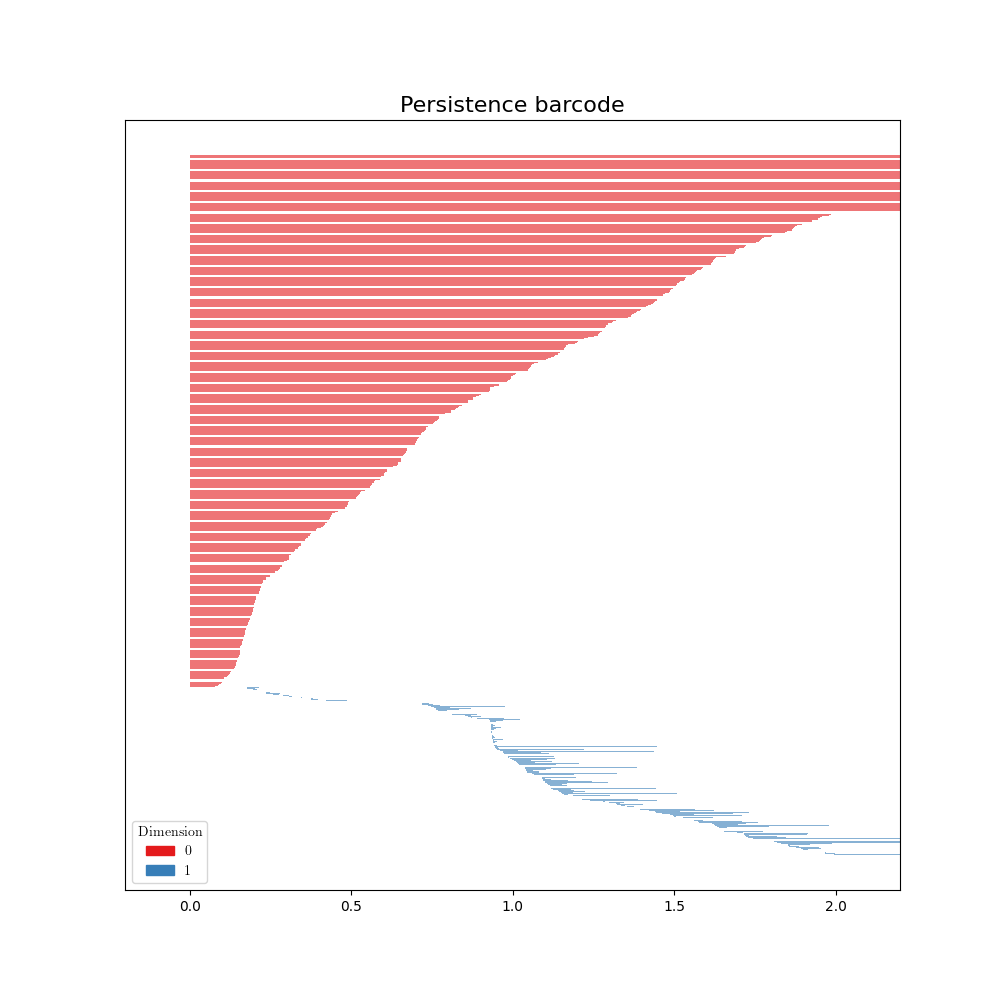}\vspace{-11pt}
                \caption{Non-contiguous US (with dates)}\label{fig:outliers_bars_dates}
	\end{subfigure}
	~
	\begin{subfigure}[b]{0.49\textwidth}
		\includegraphics[width=\textwidth]{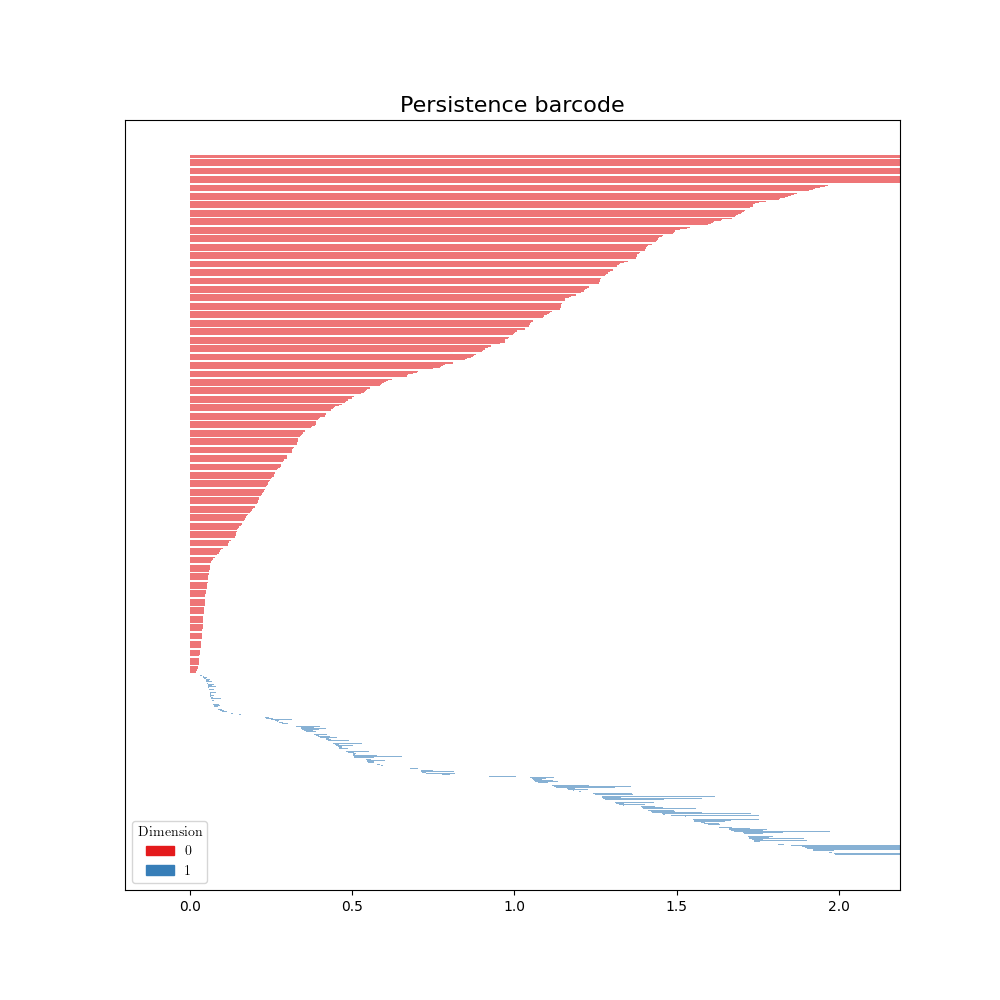}\vspace{-11pt}
                \caption{Non-Contiguous (without dates)}\label{fig:outliers_bars_nodates}
	\end{subfigure}
    \caption{The Non-Contiguous barcodes, where the $x$-axis is our time.}
    \label{fig:outliers_bars}
\end{figure}

Now looking at the features of the dataset of the Non-contiguous US, first we note that there are much fewer entries in this dataset than other datasets. Here we are only looking at 3 states/territories, so this is why we have such fewer data. We again see the increase in $S^1$ components at $t\approx 1$, which we've see before and know that this comes from the weeks and years. There are fewer $S^1$ components, than in other datasets that we have seen. This lack of $S^1$'s in due to the few entries in our dataset. We also see that these $b_1$ bars do not persist for very long in figure \ref{fig:outliers_bars_nodates}. We do see a few bars in figure \ref{fig:outliers_bars_dates}, but when we look at higher dimensional barcodes, we do not see a higher dimensional void. Therefore, we say that these bars are not coming from a higher dimensional shape.

We see that the $b_0$ bars look a little different from these two barcodes. This is coming from the dates being separated by 1 with the weeks and years. We are seeing this effect the $b_0$ components in this smaller dataset. In figure \ref{fig:outliers_bars_nodates}, we see there is a bump at $t\approx 1$. This is coming from two or more clusters joining together. Which would also explain the increase in $S^1$'s after this time. In this dataset, the shape looks like a point cloud made of a few clusters that are like arms that stick out. These arms that stick out would make the $S^1$'s.

\section{Mapper Analysis}\label{sec:Mapper}

Recall from Section \ref{sec:methods and data} that MAPPER works by first performing a principal component analysis of our dataset and determining which directions have the most variance. We have a choice of 1, 2, or 3 dimensions when we are projecting down. For the 1-dimensional option, you will not see much happening. This would just give you a line with points highlighted. For the 2-(3-)dimensional options, when we have found the two (three) directions that have the most variance, we project onto these directions using a filter function. A \textit{filter function} is a nice projection map that doesn't necessarily project onto the canonical basis directions. Once we have projected our data onto our two (three) dimensional subspace, we construct an open cover of the projected space using open balls. Then we look at the pre-image of each of these open balls, and we look at the clusters in these pre-images. If a given pre-image has a cluster, we represent it as a vertex and if two neighboring pre-images have points in common, then we attach these vertices with an edge. For our initial parameters, we set 20 consecutive intervals with 0.3 overlap. We cross referenced this with the distances used in the Vietoris-Rips construction to ensure the notion of \textit{distance} was cohesive between the two algorithms, as described in detail in Section \ref{sec:methods and data}. Below in Figure \ref{fig:circle} is an example from \cite{Mapper} to illustrate how MAPPER works.

\begin{figure}[hbt!]
    \centering
    \includegraphics[trim= 100 125 100 125,clip, scale = 0.4]{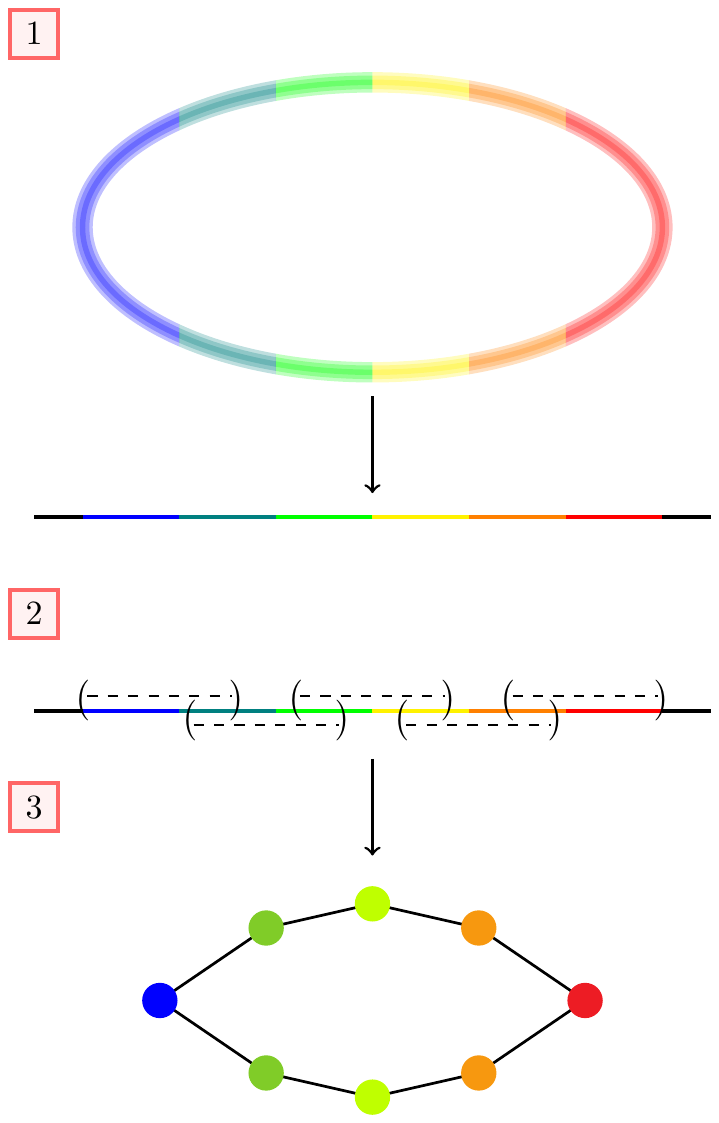}
    \caption{This is an example of a dataset that looks like a circle. We follow the steps for MAPPER: project down, construct an open cover, and then look at the components of the pre-image.}
    \label{fig:circle}
\end{figure}

For our dataset, we began by partitioning the entries by region, as defined in Section \ref{sec:methods and data}. The data were organized geographically from south to north, with entries ordered accordingly from top to bottom. For example, in the “South” dataset, Texas occupies rows 1–193, Florida occupies rows 194–387, and Virginia occupies rows 2522–2715.

Our analysis proceeds by region, with particular attention to the mapper clustering results. The interpretation relies heavily on the figures, especially the node colorings. In mapper, node colors are assigned automatically based on the row indices of the data. The lightest color (yellow) corresponds to rows with the lowest (or highest) indices, and the darkest color (dark purple) corresponds to rows with the highest (or lowest) indices. Intermediate colors, such as blue or green, may indicate that all elements in a node come from a single contiguous block of rows, or that a node contains both low- and high-index rows. In each subsection, we specify precisely which rows correspond to which colors.

We applied the clustering algorithm to both the two-dimensional and three-dimensional projections of the dataset. In a separate analysis, we removed the temporal variables (year and week number), referring to this modified dataset as the “without dates” version. We then compared and contrasted the results of the “without dates” dataset with those of the complete dataset for both the two-dimensional and three-dimensional projections. This second analysis was conducted to determine whether the removal of date information affected the clustering structure or revealed additional patterns in the data.

Additionally, Table \ref{tab:mapperClusters} lists the columns selected for projection into two-dimensional and three-dimensional space, respectively. Columns 0 and 1 correspond to the year and week number, respectively. Notably, Column 2 (total death count) and Column 16 (COVID-19 death toll) emerged as particularly impactful variables. For the dataset representing the entire United States, Column 13 (deaths due to heart disease) was also prominent (see Table \ref{tab:mapperClusters}).

\begin{table}[hbt!]
    \centering
    \begin{tabular}{|l|l|l|}
        \hline
        Region/Group & 2-dim columns & 3-dim columns\\ \hline 
        West & Columns 2 and 16& Columns 2, 16, and 1 \\[10pt] \hline
        Midwest & Columns 2 and 16& Columns 2, 16, and 1 \\[10pt]  \hline
        Northeast & Columns 2 and 16& Columns 2, 16, and 1\\[10pt] \hline
        South & Columns 2 and 16& Columns 2, 16, and 1\\[10pt] \hline
        Non-Contiguous US & Columns 3 and 16& Columns 3, 16, and 0\\[10pt]
        \hline
        Entire US & Columns 13 and 10& Columns 13, 10, and 0\\[10pt]
        \hline
    \end{tabular}
    \caption{Mapper clustering details. Column 16 contained deaths due to Covid-19.}
    \label{tab:mapperClusters}
\end{table}

We repeated the analysis on the datasets with the date columns removed (see Table \ref{tab:mapperClusters_noDates}). For each region, the selected columns in the two-dimensional projection remained identical to those from the complete dataset. However, differences emerged in the three-dimensional projection, as the week column—used in the initial clustering—was no longer available. In its absence, the next most influential variable was Column 10, corresponding to deaths from other respiratory diseases.

\begin{table}[hbt!]
    \centering
    \begin{tabular}{|l|l|l|}
        \hline
        Region/Group & 2-dim columns & 3-dim columns\\ \hline 
        West & Columns 0 and 14& Columns 0, 14, and 10 \\[10pt] \hline
        Midwest & Columns 0 and 14& Columns 0, 14, and 10 \\[10pt]  \hline
        Northeast & Columns 0 and 14& Columns 0, 14, and 10\\[10pt] \hline
        South & Columns 0 and 14& Columns 0, 14, and 10\\[10pt] \hline
        Non-Contiguous US & Columns 1 and 14& Columns 1, 14, and 10\\[10pt]
        \hline
        Entire US & Columns 11 and 14& Columns 11, 14, and 10\\[10pt]
        \hline
    \end{tabular}
    \caption{Mapper clustering details - no dates datasets}
    \label{tab:mapperClusters_noDates}
\end{table}

We now turn to the analysis of the clusters shown in the MAPPER algorithm output. For each region, we examine both two- and three-dimensional projections, with and without dates, yielding four images per region. Each region is discussed individually, with overarching conclusions presented at the end of this section. 

Note that all the figures for the MAPPER analysis are in Appendix \ref{sec:Mapper_figs}.

\subsection{Whole US}
When analyzing the mapper clustering for the entire United States (Figure \ref{fig:whole_US_data}), we observe that node colors correspond to the date. The darkest nodes (purple) represent 2020, with progressively lighter colors corresponding to later dates, culminating in yellow for December 2023. In all subfigures of Figure \ref{fig:whole_US_data}, at least one distinct cluster consists entirely of purple nodes. Notably, even when the date columns are removed, the dataset still exhibits disjoint clustering that seem to be driven by the date.

\begin{figure}[!htbp]
    \centering
    \includegraphics[width=0.5\linewidth]{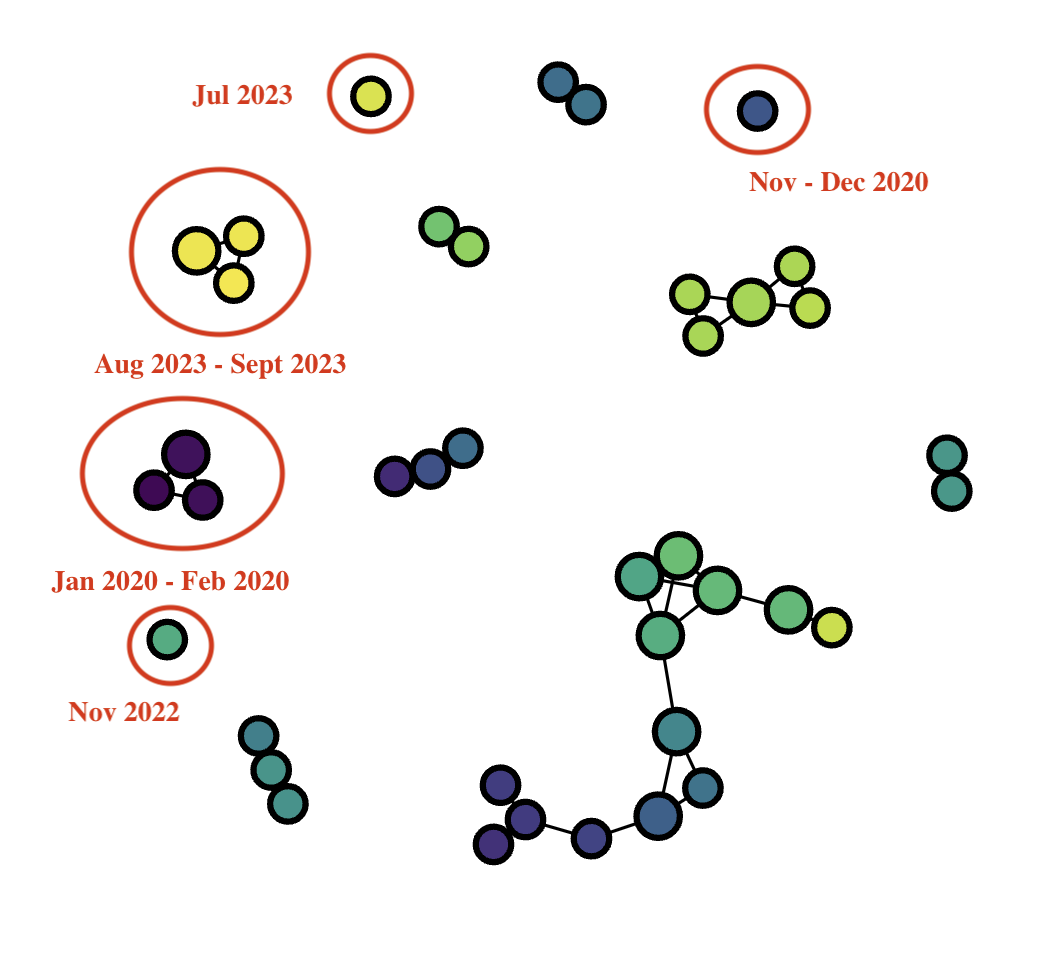}
    \caption{Whole US (without the dates)}
    \label{fig:infoWholeUS}
\end{figure}

Closer examination reveals that the first 12 weeks of 2020 consistently form either their own small cluster or, in the two-dimensional projection (Figure \ref{fig:whole_US_2D}), are connected to a green node containing rows corresponding to February and March 2023. One would assume the dates columns are the main contributing factors to the disjoint clustering that is indeed grouped by date. However, according to Table \ref{tab:mapperClusters}, the column corresponding to weeks is the third most relevant column for MAPPER. This means that in the 2-dimensional projection, the weeks column is not considered. What is even more notable, as mentioned above, is the clustering in Figures \ref{fig:infoWholeUS} and \ref{fig:whole_US_3D_noDates} is completely separated by dates, even though the dates columns were removed. This leads us to conclude an answer to Question \ref{question1} that the dates are heavily contributing to the topology of this dataset.

We wish to recall the barcode analysis in Section \ref{sec:barcode}, where we determined that this dataset looks like a point cloud with a few arms sticking out, similar to the non-contiguous states. The MAPPER analysis reflects this conclusion. In context, mortality rates across the country remained relatively consistent during the pandemic, with notable outliers occurring in the months preceding and following the pandemic period.

\subsection{West}\label{sec:west_mapper}

 The three-dimensional clustering of the West is particularly noteworthy, as all outlying points, including the peripheral clusters, are colored deep purple. This pattern suggests that either a single state or a specific year behaves very differently from the rest of the dataset. The two-dimensional mapper output exhibits a similar pattern.

\begin{figure}[!htbp]
    \centering
    \includegraphics[width=0.5\linewidth]{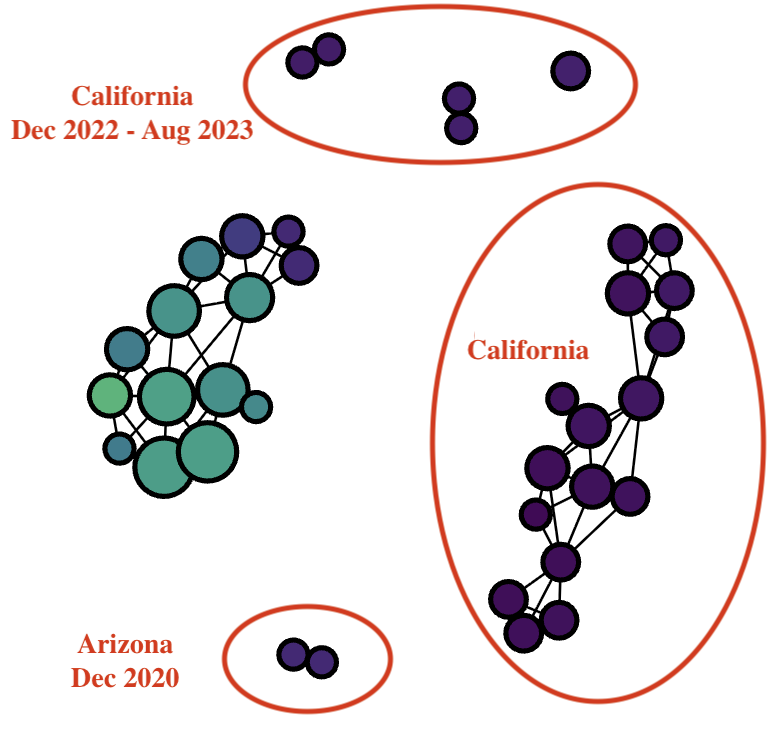}
    \caption{The West 2D (without the dates)}
    \label{fig:infowest}
\end{figure}

Closer examination of the West’s mapper clustering diagrams reveals that the state of California accounts for nearly the entire outlying set. In Figure \ref{fig:west_data}, the nodes corresponding to California are uniformly dark purple, with the exception of two nodes connected by a single edge. In all four panels of Figure \ref{fig:west_data}, these two nodes correspond to rows 246–249, representing the state of Arizona during November and December 2020. This outlying set is most prominent in Figure \ref{fig:infowest}, where the relevant node appears slightly blurred. Thus, we may pose an answer to Question \ref{question1} as the state of California is heavily contributing to the shape of this dataset.

When we looked at the west in the barcode analysis, we determined that all the data points are close together and that there are circles that do not persist for long. This led us to conclude that our dataset looks like a collection of small circles that are all close together. This resembles something like a radio tower. We can also see this in the MAPPER in figure \ref{fig:west_data}. We can see that many of the points are close together, and there are few outliers. We also see in the MAPPER that the datasets look kind of like radio towers, especially in figure \ref{fig:west_3D}. In the barcode for the west, we saw that all of the data points had merged into one connected component fairly quickly. This was the major feature that told us all our points were close together. The data points that cause the trailing off in this dataset are the data points coming mostly from California. California's data is so separate from the other data points that it had its own clusters. We again see this in the MAPPER analysis in figure \ref{fig:west_2D}. The data joining up so quickly and having California being so far removed tells us that during the pandemic mortality in the other states in the west were all pretty close to each other and the amount/pattern of California's mortality was much different than the rest of the states. Therefore, we conclude that the shape of the data makes sense and the main body of the data is formed from all the states except California. Furthermore, California's data is removed from the main body of the data and forms its own clusters.

\subsection{Midwest}

For the Midwest, Figure \ref{fig:midwest_data}, the yellow and very light green nodes correspond to North Dakota and Minnesota, respectively. Unlike the West and South, no single state emerges as a clear outlier. However, the clustering suggests a notable relationship between data from Minnesota and North Dakota. Rows corresponding to Illinois frequently appear in outlying nodes, while Missouri and Indiana often occur together. Additionally, South Dakota, North Dakota, Minnesota, and Wisconsin tend to be grouped within the same nodes. A preliminary examination indicates that these four states rarely share small nodes with the other Midwest states. A potential explanation for this behavior is geographic proximity, which could lend towards a partial answer to Question \ref{question1} for the Midwest dataset.

Recall from the barcode analysis, we determined that the midwest dataset does not depend on the dates too much. This is significant because this tells us that the true shape of the dataset still shows up in the higher dimensions produced by the weeks and years. We concluded that the shape of the dataset is something like a point cloud with few clusters spread out in a linear way away from a `center'. This is similar to what we see in the MAPPER figures, especially in the figures of Figure \ref{fig:midwest_data}. We can see a few clusters that stick out from a center. There are few outliers, which also contribute to trailing in the barcodes, coming from many different states in different choices of 2D or 3D, or with/without the dates. Putting all of this information in context, we arrive at the conclusion that the Midwest region, like the West region, experience few instances of mortality being so far removed from a main body of data points. The few instances of clusters away from the main body of data points come from data at January 2020 and August/September 2023. This is telling us that mortality during the pandemic was much different from before the pandemic started and from when the cases started to slow at the end of 2023. Therefore, we say that the midwest had mortality that, while the pandemic was going on, stayed close to each other. And when we get to August/September 2023, mortality reached a `new normal' that was different from the mortality before the pandemic began.

\begin{figure}[!htbp]
    \centering
    \includegraphics[width=0.5\linewidth]{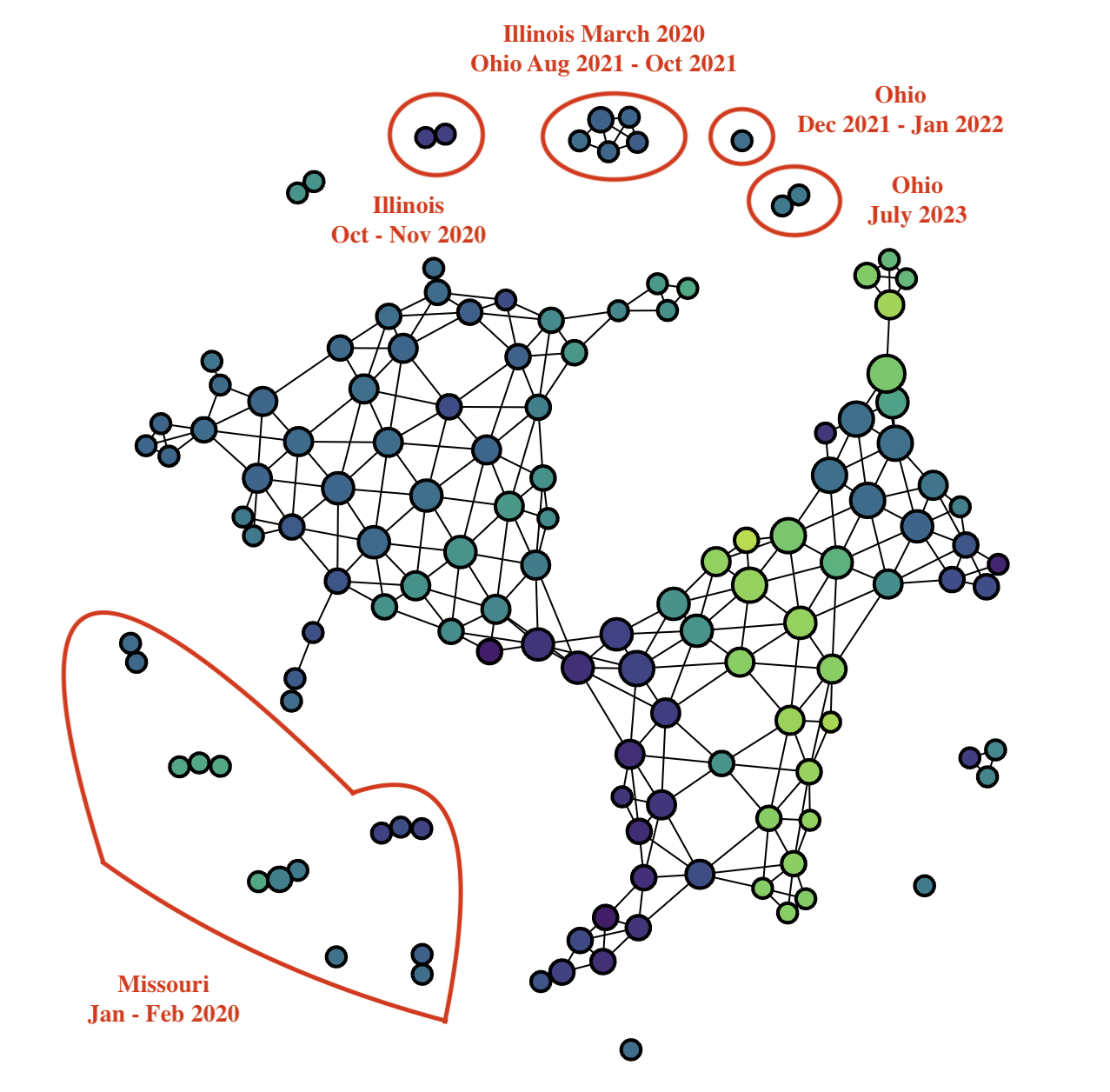}
    \caption{The Midwest (without the dates)}
    \label{fig:infomidwest}
\end{figure}

\subsection{Northeast}
\begin{figure}
    \centering
    \includegraphics[width=0.5\linewidth]{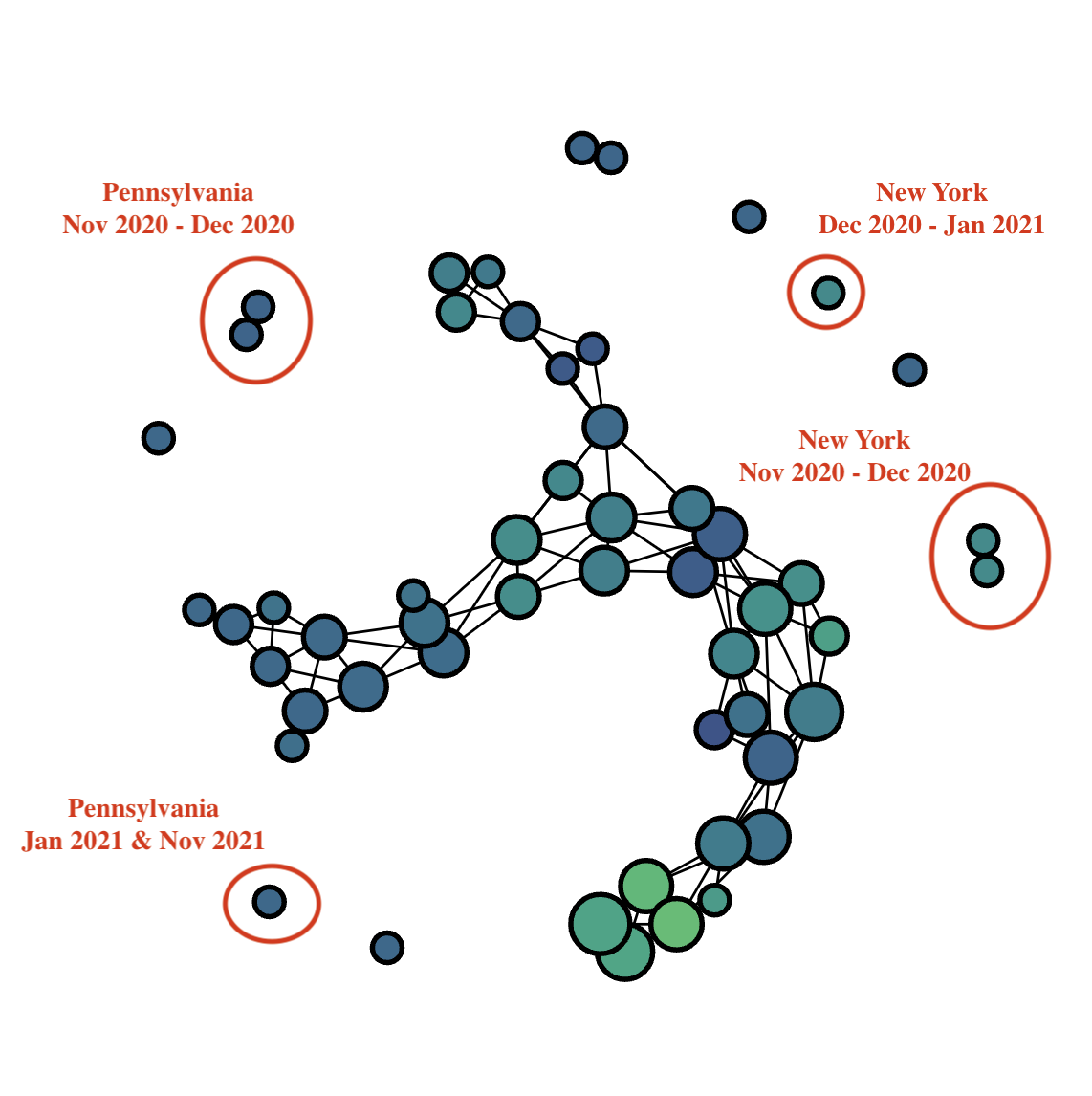}
    \caption{Northeast (without the dates)}\label{fig:infnortheast}
\end{figure}

The yellow and light green nodes correspond to Maryland, Washington, D.C., and Delaware. Notably, the yellow nodes, which represent Maryland, are never identified as outliers. In the 2D projection for the dataset without the dates (Figure \ref{fig:infnortheast}), yellow nodes are absent altogether. This suggests a high degree of integration or similarity between Maryland and other states in the dataset. The darkest purple nodes, which represent Vermont, also do not appear as outliers. Most nodes display a blend of colors, indicating that they are composed of data from multiple states. This lack of distinct outliers from individual states may be attributed to the relatively small size and close geographic proximity of the Northeastern states, which could contribute to their data being more interwoven.

Upon closer analysis of the outlying clusters, we see that New York and Pennsylvania during the winter months consistently forms an outlying node or cluster in all four cases in Figure \ref{fig:northeast_data}. The unlabeled outlying clusters in Figure \ref{fig:infnortheast} are various months for the states of Pennsylvania and New York. We also see New York City during the winter months forming a distinct cluster in all cases except 2D without dates. This reflects the same conclusion from the barcode analysis done for the Northeast from Section \ref{sec:barcode}.

In the barcode analysis, we noted that there was a separation of our clusters in the dataset that includes dates because of the weeks and years. We concluded that since this was not the case in the dataset without the dates, that there was probably some clusters that were due to spikes during particular times of the year that repeated every year around the same week number. We see this in the MAPPER with the outlying clusters of New York, Pennsylvania, and New York City in the winter time. Therefore, our conclusion that the dataset resembles a point cloud with some striations is consistent to what we see in the MAPPER analysis. These observations make sense when we put this information into context. During the COVID-19 pandemic New York City was an epicenter of many different things happening and therefore it makes sense that it would be an outlier in our analyses.

\subsection{South}
For the South, Figure \ref{fig:south_data}, the dark purple nodes correspond to Texas and Florida. In each of the four cases, these two states consistently form their own distinct cluster. This separation is particularly notable given their lack of geographic proximity. However, their alignment can be explained by the fact that Texas and Florida adopted nearly identical COVID-19 lockdown procedures, which set them apart from the rest of the southern states.

From the barcode analysis, we noted that the south dataset is structured differently from the rest of the region's datasets. We can see this difference in the barcodes (Figure \ref{fig:South_bars}) when compared to any of the other datasets. The tells us that the south treated the pandemic much different than other regions. In the barcode analysis, we concluded that the dataset looks like a point cloud with few, if any clusters. But in the MAPPER analysis, this is not what we see. We can clearly see two distinct clusters happening, see Figure \ref{fig:south_data}. So, how do we reach two different conclusions for this dataset? We reach this conclusion because the second smaller cluster is still away from the main cluster in a normal way. The second cluster is coming from Texas and Florida. When our filtration goes through these data points we saw that all the points were spaced out in a normal way. This means that Texas and Florida still follow this normal distribution while also being far enough removed from the main body to form their own cluster before joining with the main body of points. All of this information tells us that mortality in the south grew steadily at the beginning of the pandemic and decreased steadily towards the end of the pandemic forming a main cluster with few outliers. The `spikes' in mortality were pretty much concentrated in Texas and Florida which caused them to be so far removed from the main body of points, but to still be close enough to join with each other.

\begin{figure}[!htbp]
    \centering
    \includegraphics[width=0.5\linewidth]{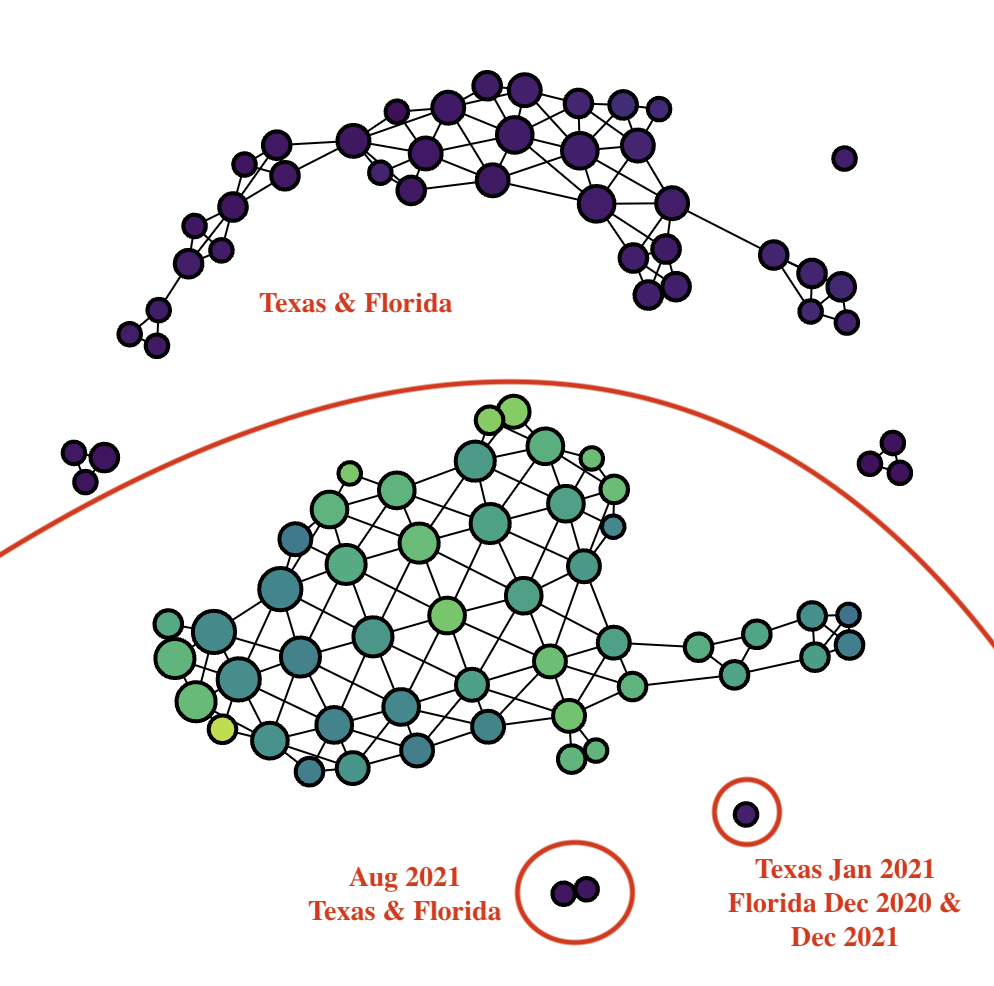}
    \caption{The South (without the dates)}
    \label{fig:infosouth}
\end{figure}

\subsection{Non-Contiguous US}
The outlier group consists of Alaska, Hawaii, and Puerto Rico. Alaska occupies rows 0–195, appearing as the darkest purple in the 2020 data and a lighter purple in the 2023 data. Hawaii occupies rows 196–389, corresponding to darker blue in 2020 and lighter blue in 2023. Puerto Rico occupies rows 390–583, with light green representing 2020 and yellow representing 2023. As shown in Figure \ref{fig:outliers_data}, the yellow and light green nodes consistently form clusters disjoint from the blue and purple nodes. Given that these regions do not share borders and are geographically distant, it is unsurprising that the projections produce highly disconnected clusters. Indeed, we observe separate clusters corresponding to each state or territory, as expected.

When we did the barcode analysis, we found that our dataset looks like a point cloud with a few arms sticking out to form the $S^1$ components that we saw in Figure \ref{fig:outliers_bars}. In the MAPPER analysis we see similar information, especially in Figure \ref{fig:outliers_data} except for Figure \ref{fig:outliers_3D}. Here we see many disjoint clusters, but when we go down one more dimensions (Figure \ref{fig:outliers_2D}) we see that these disjoint clusters are now apart of a cluster like we describe in the barcode analysis. This is telling us that the column 0 (the year) is holding a lot of information. It's holding enough to separate clusters like it has in Figure \ref{fig:outliers_3D} (see Table \ref{tab:mapperClusters}). Now to put all of this in context we see that the dataset for the non-contiguous states and territories looks like a point cloud with a few arms sticking out. Each of these arms are coming from the states' differences. We conclude that this means during the pandemic our states' mortality stayed consistent from state to state in these states.

\begin{figure}[!htbp]
    \centering
    \includegraphics[width=0.5\linewidth]{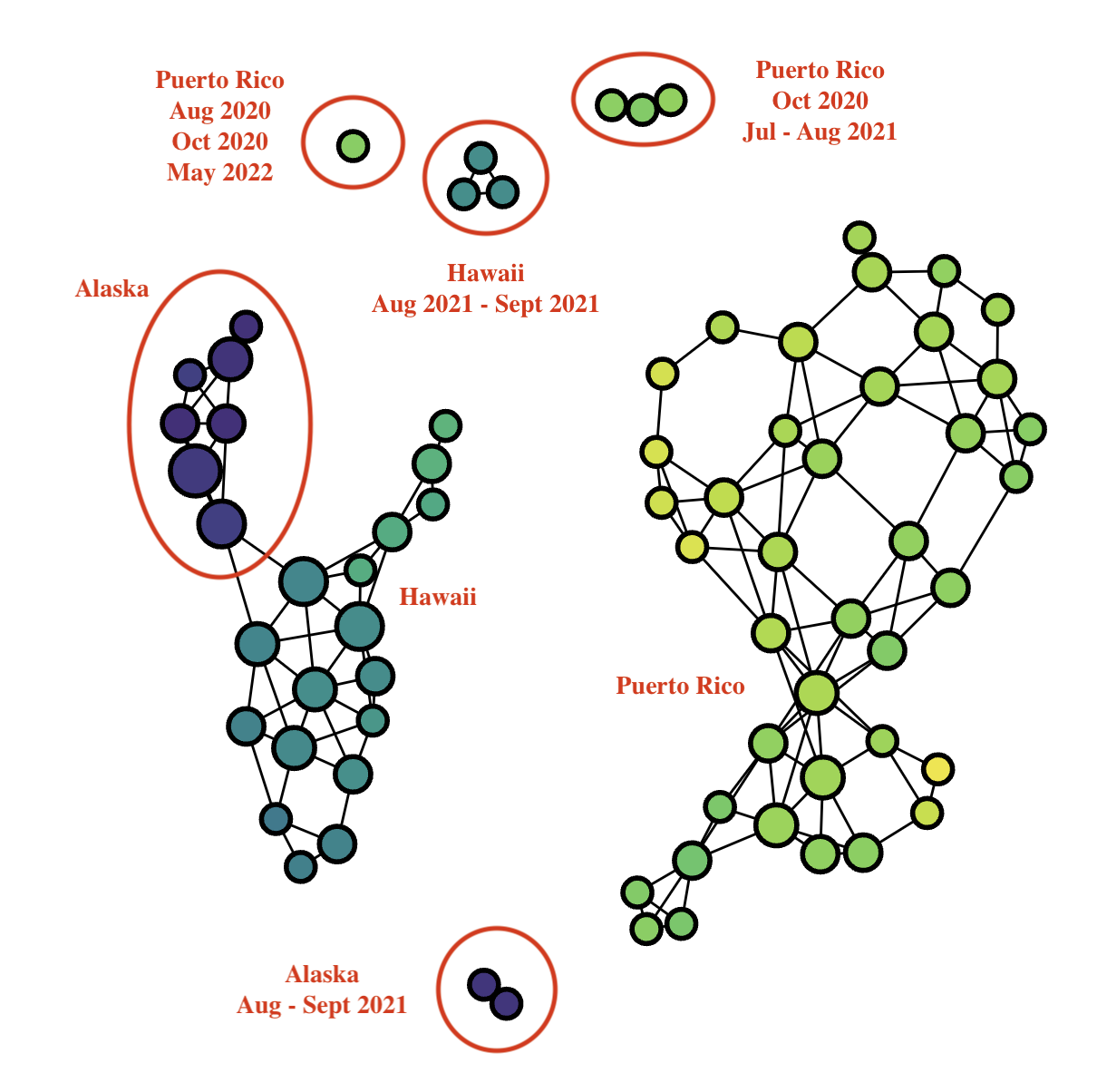}
    \caption{The Non-Contiguous US (without the dates)}
    \label{fig:infononcontiguous}
\end{figure}

\section{Conclusions}\label{sec:conclusions}

The main goal of this paper was to answer questions about contributing factors to the topology of our dataset using barcode and MAPPER analysis, respectively. We overview our answers to our main three questions, restated below for the reader, in this section. 

\newtheorem*{thm_1}{Question~\ref{question1}}
\begin{thm_1}  What are the main contributing factors to the overall shape of the dataset for the United States versus the geographical regions?   \end{thm_1}

Our original assumption was the dates or geographical proximity would be the contributing factors to the topology of our dataset. We found this to be the case with the Whole US, Non-Contiguous US, and the Midwest. When analyzing the barcode and MAPPER outcomes, we found that the dates only heavily contributed to the clustering of the Whole US dataset. As expected, the Non-Contiguous US was clustered by geographic region. The Midwest also gave us geographic proximity as the main contributor to the clusters of its dataset. 

An unexpected contributor to the topology of the South and West was lockdown procedures during the pandemic \cite{covid_lockdowns}. We found in the South barcodes and MAPPER analyses that COVID-19 lockdown procedures are the contributing factor for the topology of this dataset. Texas and Florida had different procedures compared to the rest of the states, which we believe is the reason they form their own distinct cluster. The West had an intriguing outlier as the state of California. We infer this is because California had the most strict and longest enforced lock down procedures, in stark contrast to Texas and Florida, see \cite{covid_lockdowns}.

The Northeast primary outliers were from Pennsylvania, New York, and New York City in the winter, which leads to us believing the population density of this region plus the harsh winters is the primary factor contributing to the topology of the dataset. 

\newtheorem*{thm_2}{Question~\ref{question2}}
\begin{thm_2}  Are the topological features depicted in the barcode analysis also depicted in the MAPPER projection?  \end{thm_2}

Yes, we find this to be true in all datasets. In all the datasets that we looked at in this paper, we found that the barcode analysis and the MAPPER analysis both summarize the datasets. We saw in the west that our radio tower structure from the barcode analysis could be clearly seen in the MAPPER projection. We inferred from the Non-Contiguous US, Northeast, Whole US, and Midwest barcode analyses that our dataset looks like a point cloud with varying degrees of interesting topological information the MAPPER analyses reflected these observations. We also saw with the South dataset that the barcode and the MAPPER had similar analyses: a normal distribution of points in both the clusters in the MAPPER.

\newtheorem*{thm_3}{Question~\ref{question3}}
\begin{thm_3}
    How large of an impact did the COVID-19 mortality rate have on the overall shape of the various datasets?
\end{thm_3}

The COVID-19 mortality rate had the largest impact on the Whole US dataset. This was clear, as even when we examined the dataset without dates, we still saw three definitive clusters as before, during, and after the pandemic. In the datasets corresponding to geographic regions, we notice the mortality rate had an effect, but not a big enough impact to cause clustering and disconnected components. 

The COVID-19 mortality rate impacted the West and South through the lens of lockdown procedures, while it affected the Northeast more severely during the winter. However, the topology of the Midwest dataset seemed to rely on geography rather than seasonal or otherwise. As expected, the Non-Contiguous US was completely separated by state. 

In conclusion, our exploration of the dataset provided by the CDC recording mortality during the pandemic gave us various results based on geographic region. We have found that the dataset gives us insight into how these regions responded to the COVID-19 pandemic and how different geographic features told contributed to the pandemic. From our analyses, we can see that, for example, the South responded much different than the Northeast. We were able to use tools from topological data analysis to see the different ways that these regions treated the pandemic. In future work, we want to explore using more recent tools in topological data analysis to infer more information about these regions and the impact that the COVID-19 pandemic had on their mortality. We hope that the information provided and the analyses conducted using TDA on these datasets can be used as an example to help other researchers incorporate TDA in their work.

\newpage

\nocite{*}
\bibliography{Bibliography}

\bibliographystyle{plain}

\makeatletter
\@setaddresses
\let\@setaddresses\relax
\makeatother

\newpage

\appendix

\section{MAPPER figures}\label{sec:Mapper_figs}

In this appendix, we include the MAPPER projections for the different regions. We have the 2-dimensional projection with and without the dates, and the 3-dimensional projection with and without the dates.

\begin{figure}[!htbp]
	\begin{subfigure}[b]{0.4\textwidth}
		\includegraphics[width=\textwidth]{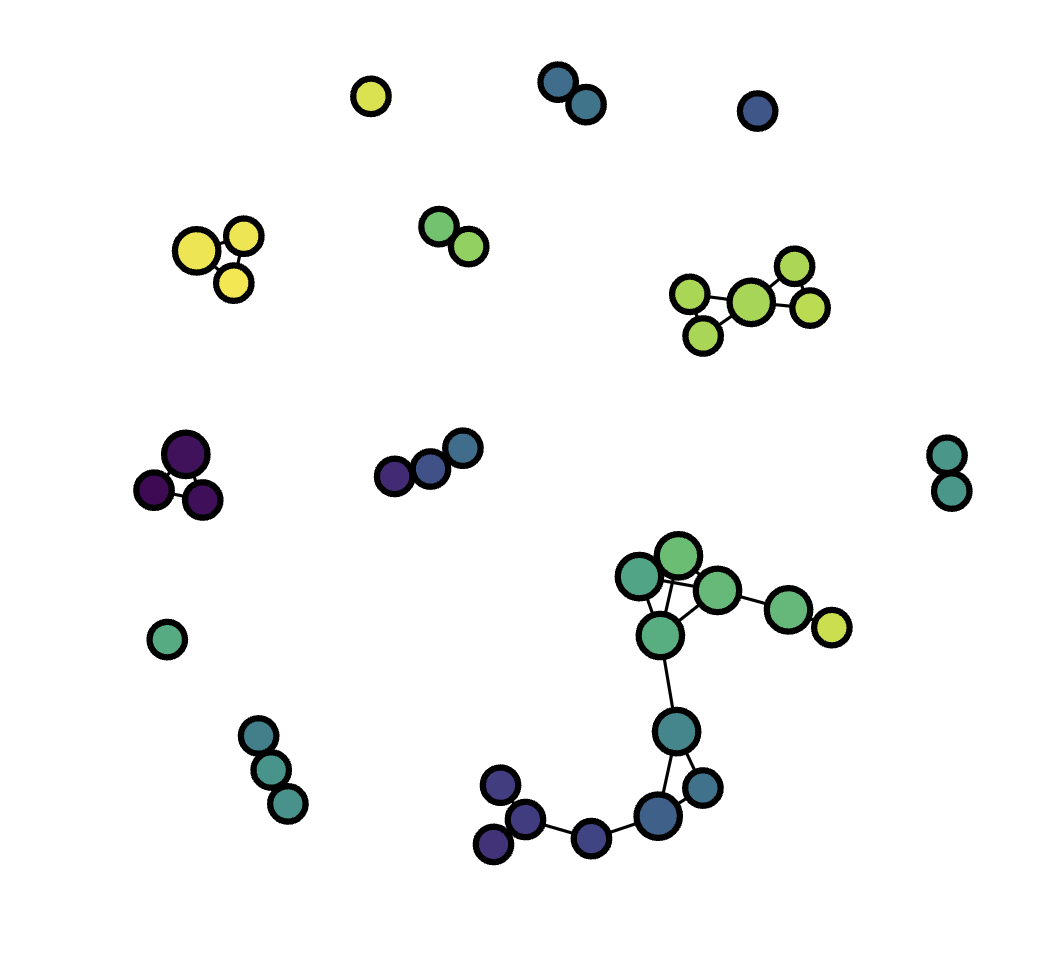}\vspace{-11pt}
                \caption{Whole US 2D (without the dates)}\label{fig:whole_US_2D_noDates}
	\end{subfigure}
	\begin{subfigure}[b]{0.35\textwidth}
		 \includegraphics[width=\textwidth]{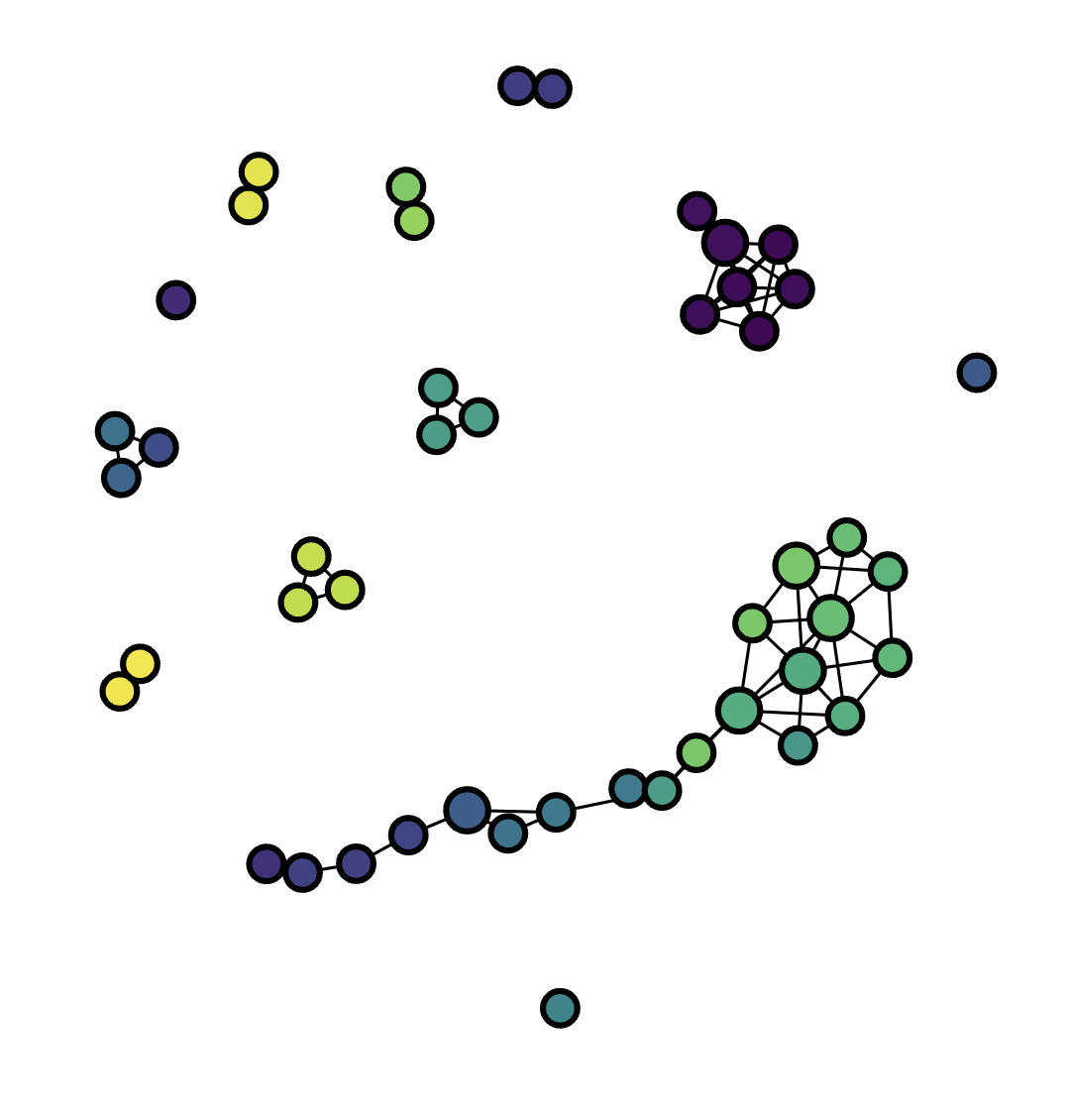}\vspace{-11pt}
                \caption{Whole US 3D (without the dates)}\label{fig:whole_US_3D_noDates}
	\end{subfigure}
    
	\vskip3mm

        \begin{subfigure}[b]{0.4\textwidth}
		\includegraphics[width=\textwidth]{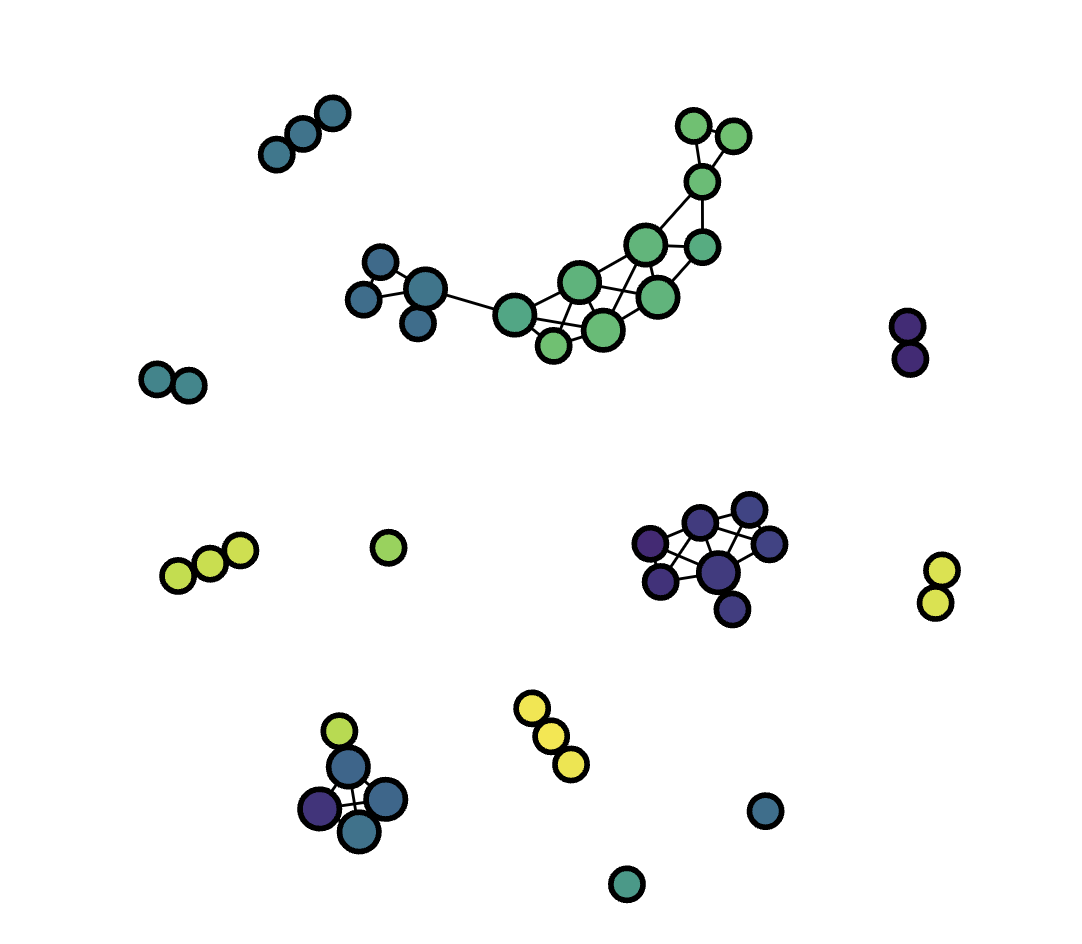}\vspace{-9pt}
                \caption{Whole US 2D (with the dates)}\label{fig:whole_US_2D}
	\end{subfigure}
	~
	\begin{subfigure}[b]{0.4\textwidth}
		\includegraphics[width=\textwidth]{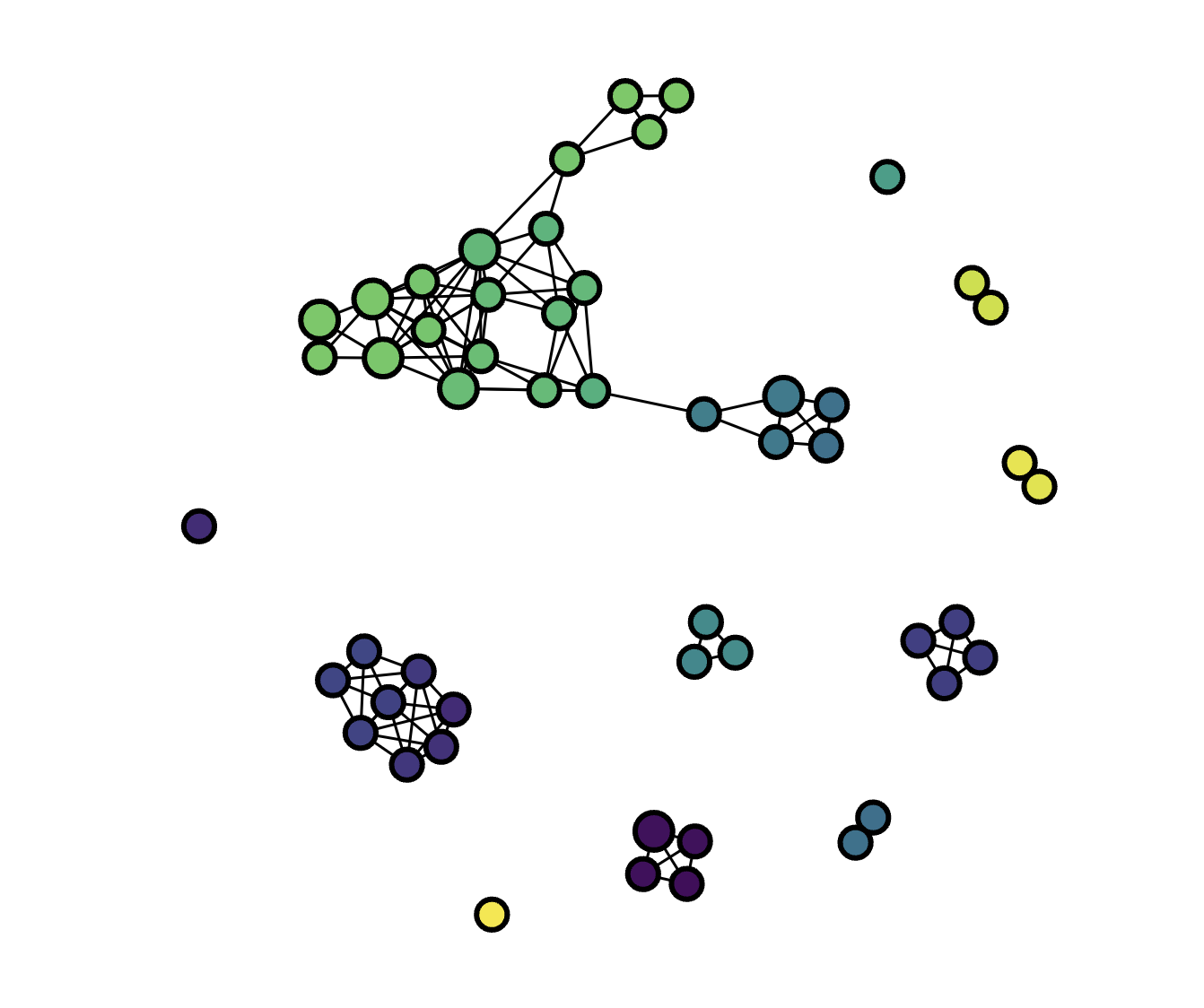}\vspace{-9pt}
                \caption{Whole US 3D (with the dates)}\label{fig:whole_US_3D}
	\end{subfigure}
        \vskip3mm
	\caption{Whole US Clustering via Mapper}\label{fig:whole_US_data}
\end{figure}

\begin{figure}[!htbp]
	\begin{subfigure}[b]{0.4\textwidth}
		\includegraphics[width=\textwidth]{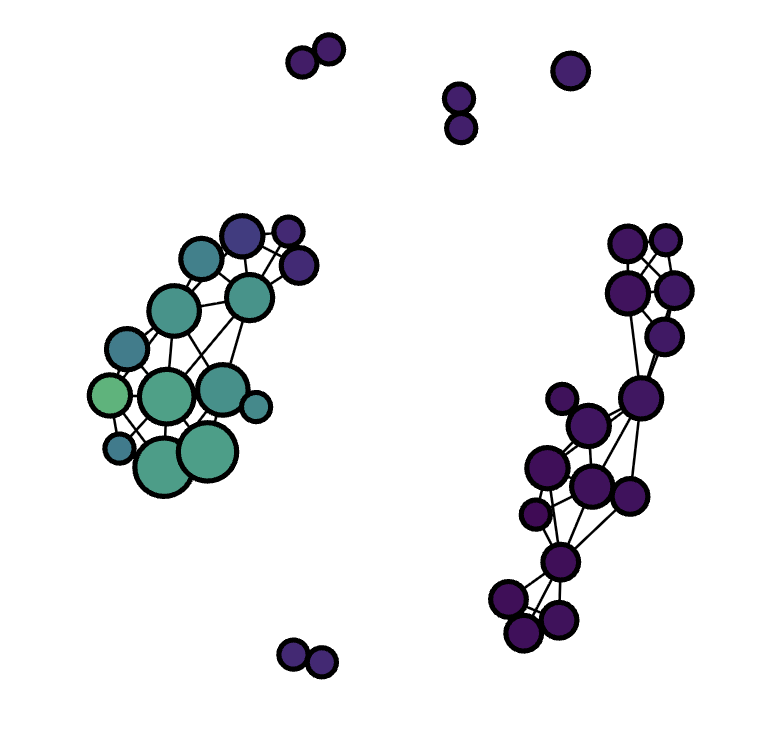}\vspace{-11pt}
                \caption{West 2D (without the dates)} \label{fig:west_2D_noDates}
	\end{subfigure}
	\begin{subfigure}[b]{0.4\textwidth}
		 \includegraphics[width=\textwidth]{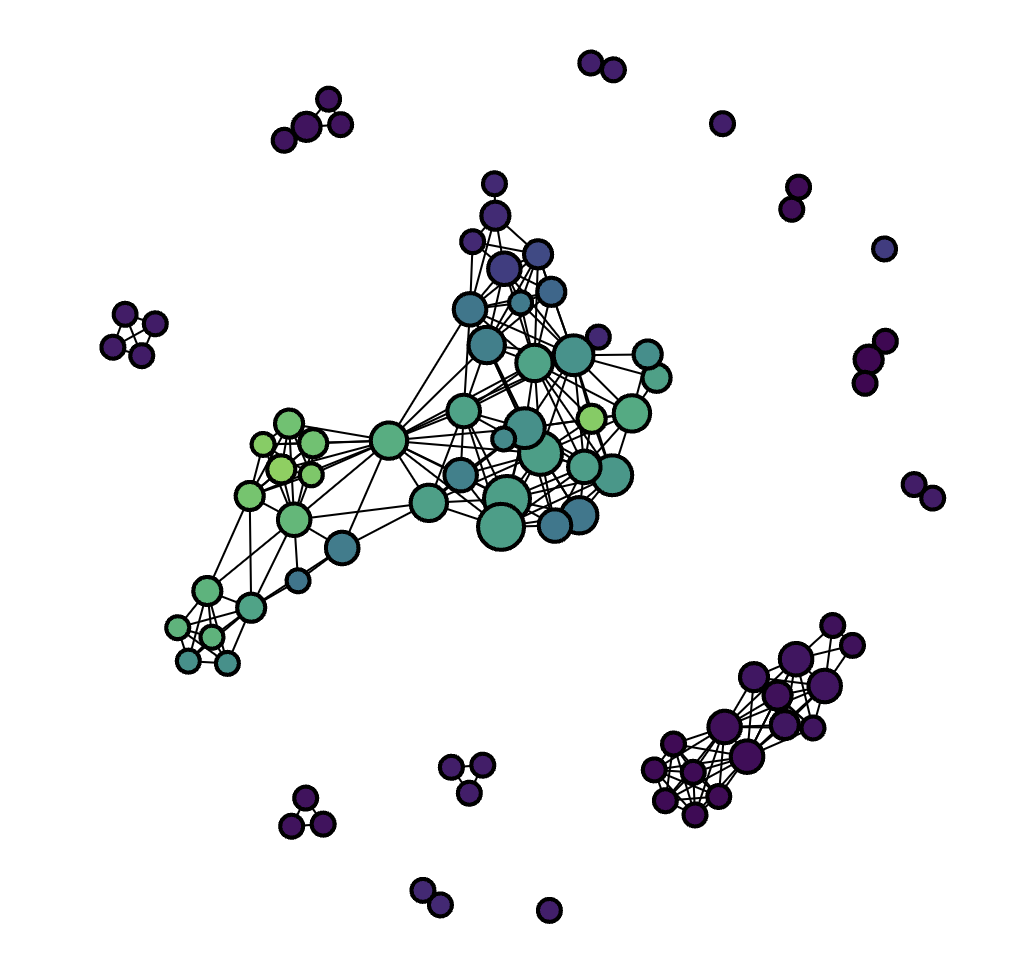}\vspace{-11pt}
                \caption{West 3D (without the dates)}\label{fig:west_3D_noDates}
	\end{subfigure}
    
	\vskip3mm

        \begin{subfigure}[b]{0.38\textwidth}
		\includegraphics[width=\textwidth]{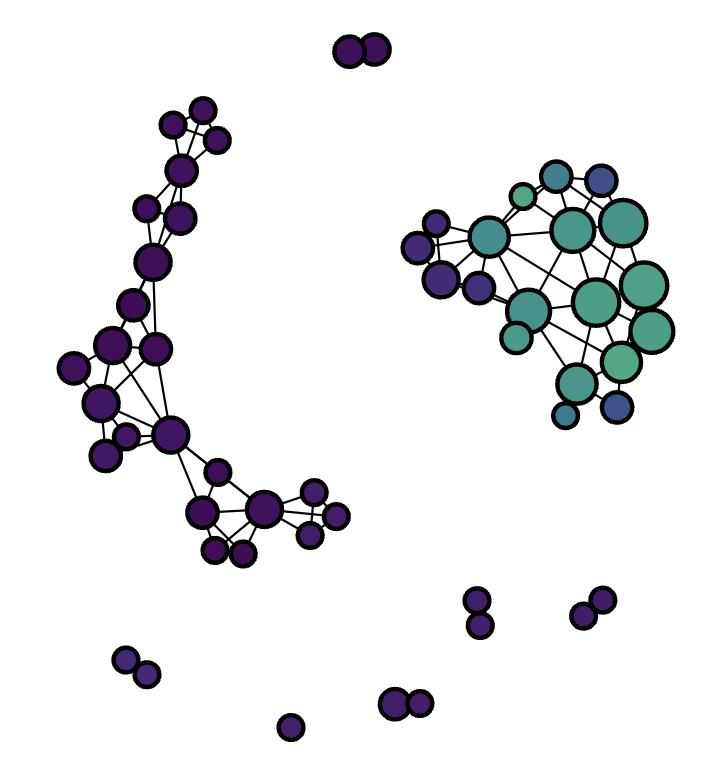}\vspace{-11pt}
                \caption{West 2D (with the dates)}\label{fig:west_2D}
	\end{subfigure}
	~
	\begin{subfigure}[b]{0.48\textwidth}
		\includegraphics[width=\textwidth]{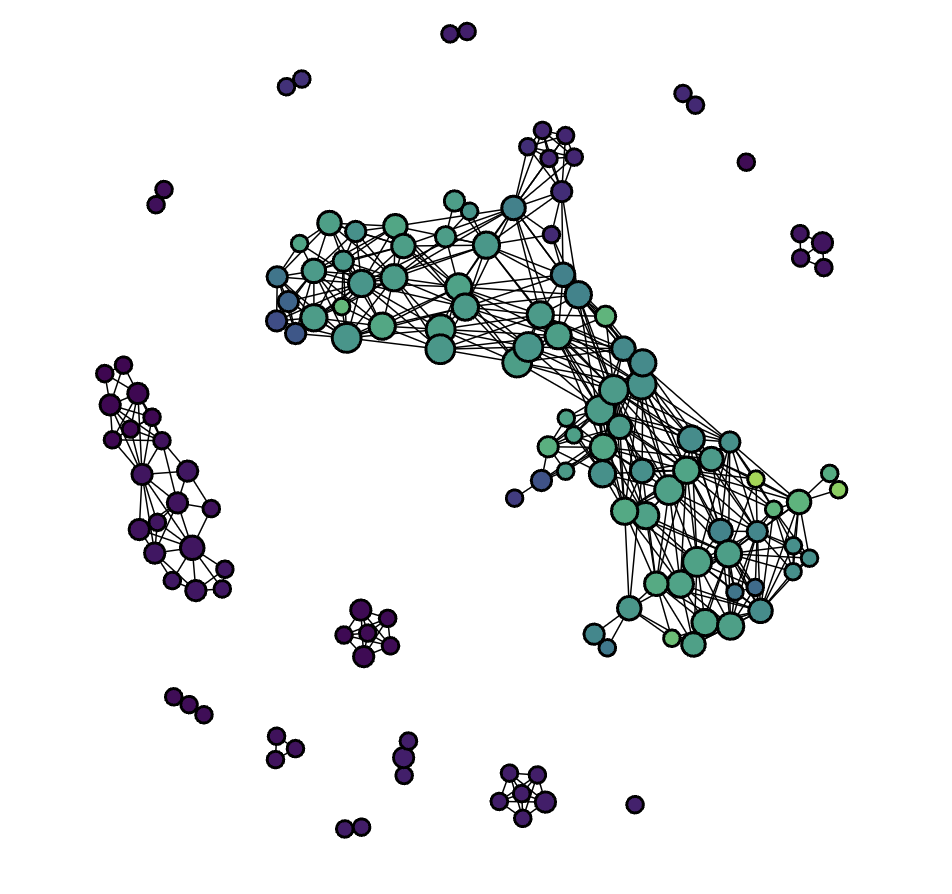}\vspace{-11pt}
                \caption{West 3D}\label{fig:west_3D}
	\end{subfigure}
        \vskip3mm
	\caption{West Clustering via Mapper (with the dates)}\label{fig:west_data}
\end{figure}

\begin{figure}[!htbp]
	\begin{subfigure}[b]{0.4\textwidth}
		\includegraphics[width=\textwidth]{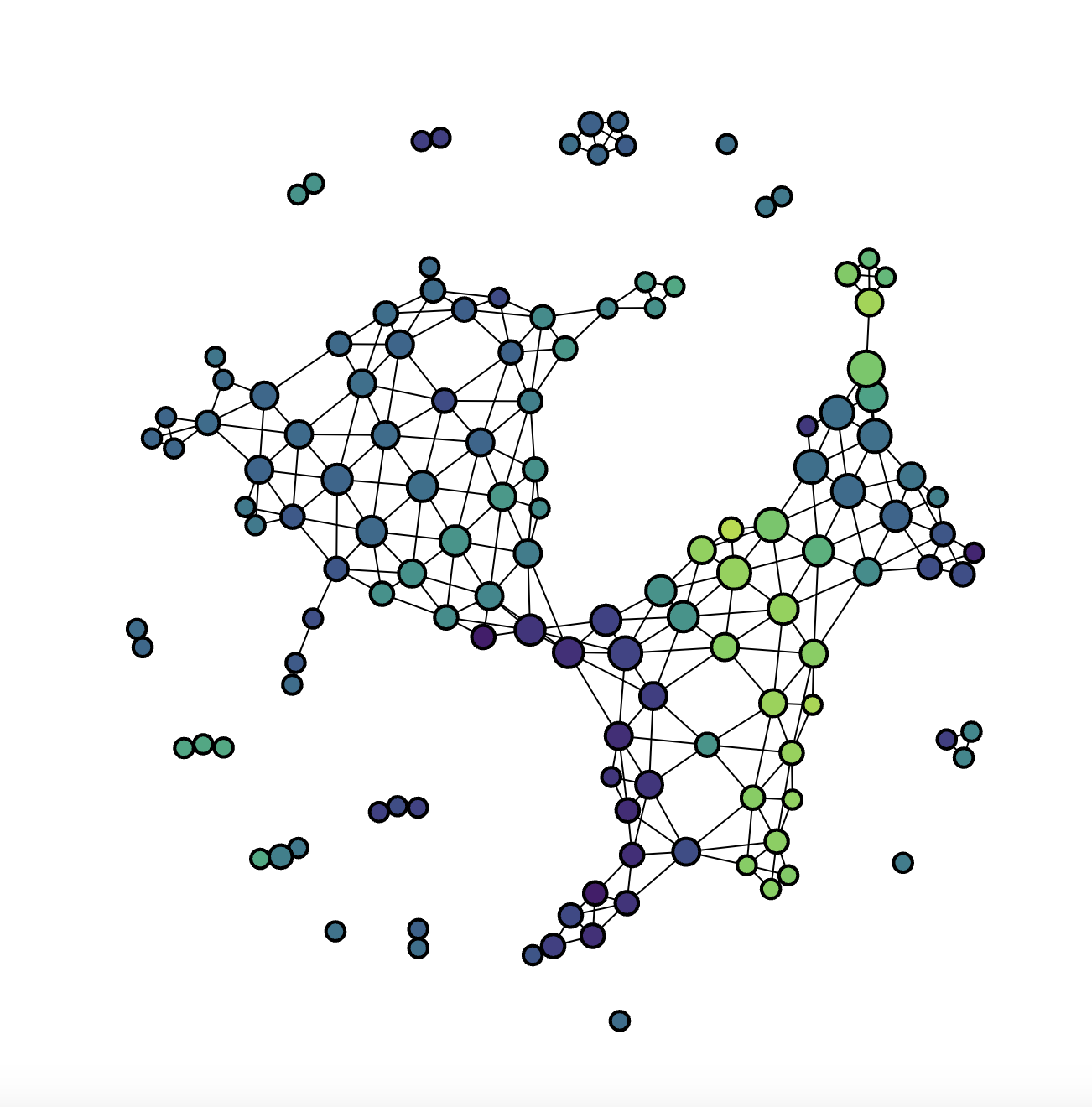}\vspace{-11pt}
                \caption{Midwest 2D (without the dates)} \label{fig:midwest_2D_noDates}
	\end{subfigure}
	\begin{subfigure}[b]{0.4\textwidth}
		 \includegraphics[width=\textwidth]{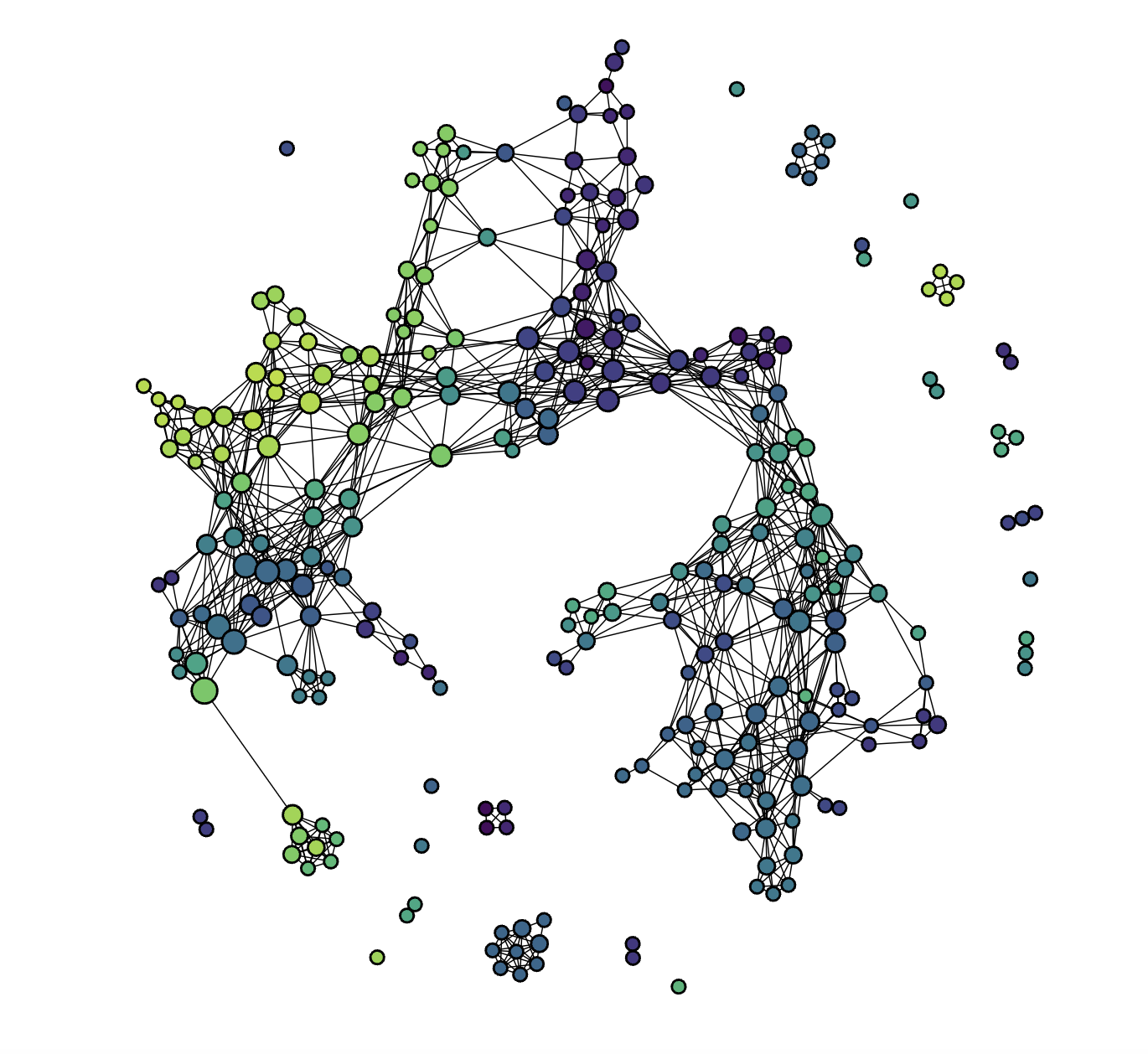}\vspace{-11pt}
                \caption{Midwest 3D (without the dates)}\label{fig:midwest_3D_noDates}
	\end{subfigure}
    
	\vskip3mm

        \begin{subfigure}[b]{0.43\textwidth}
		\includegraphics[width=\textwidth]{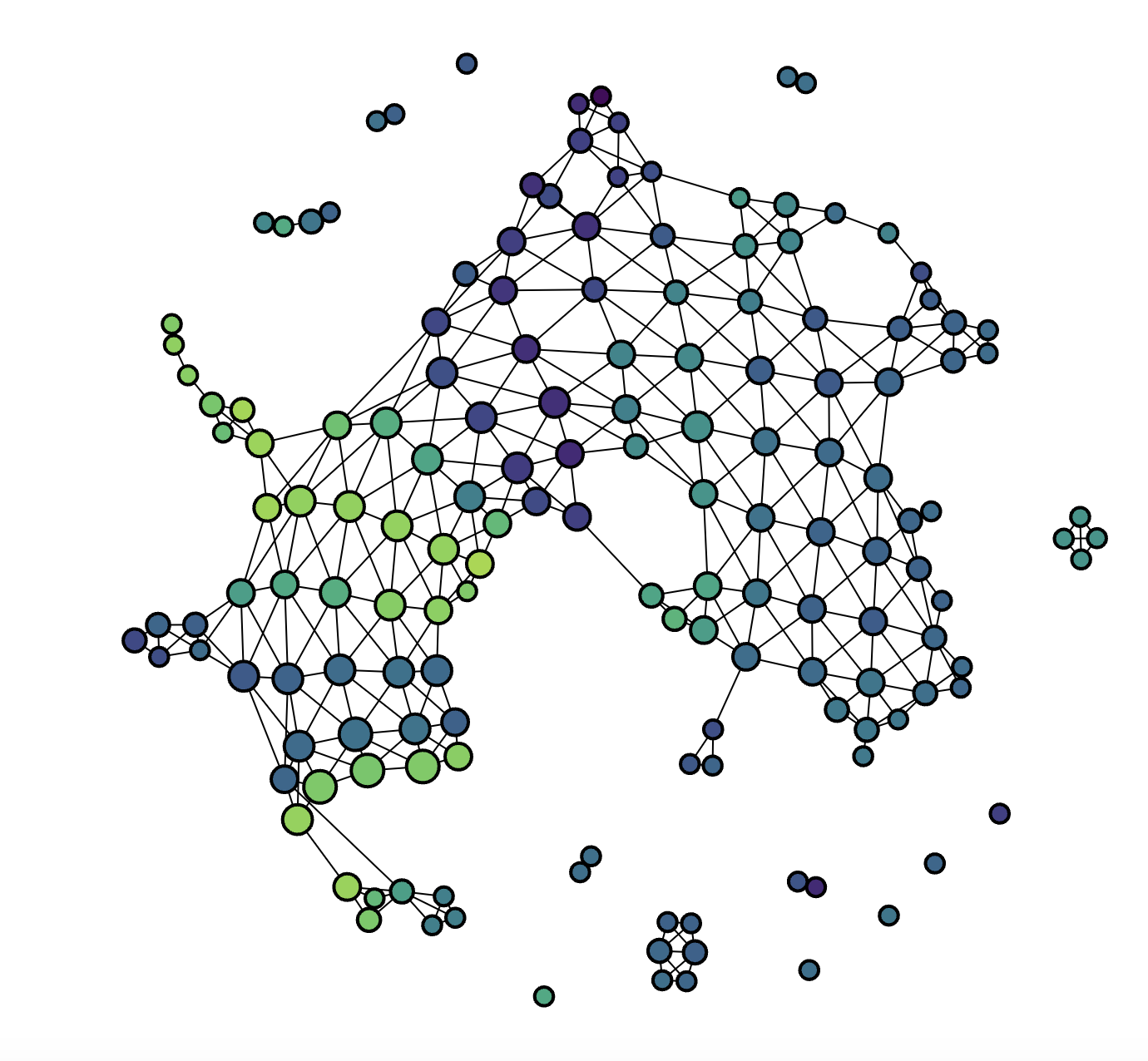}\vspace{-11pt}
                \caption{Midwest 2D (with the dates)}\label{fig:midwest_2D}
	\end{subfigure}
	~
	\begin{subfigure}[b]{0.48\textwidth}
		\includegraphics[width=\textwidth]{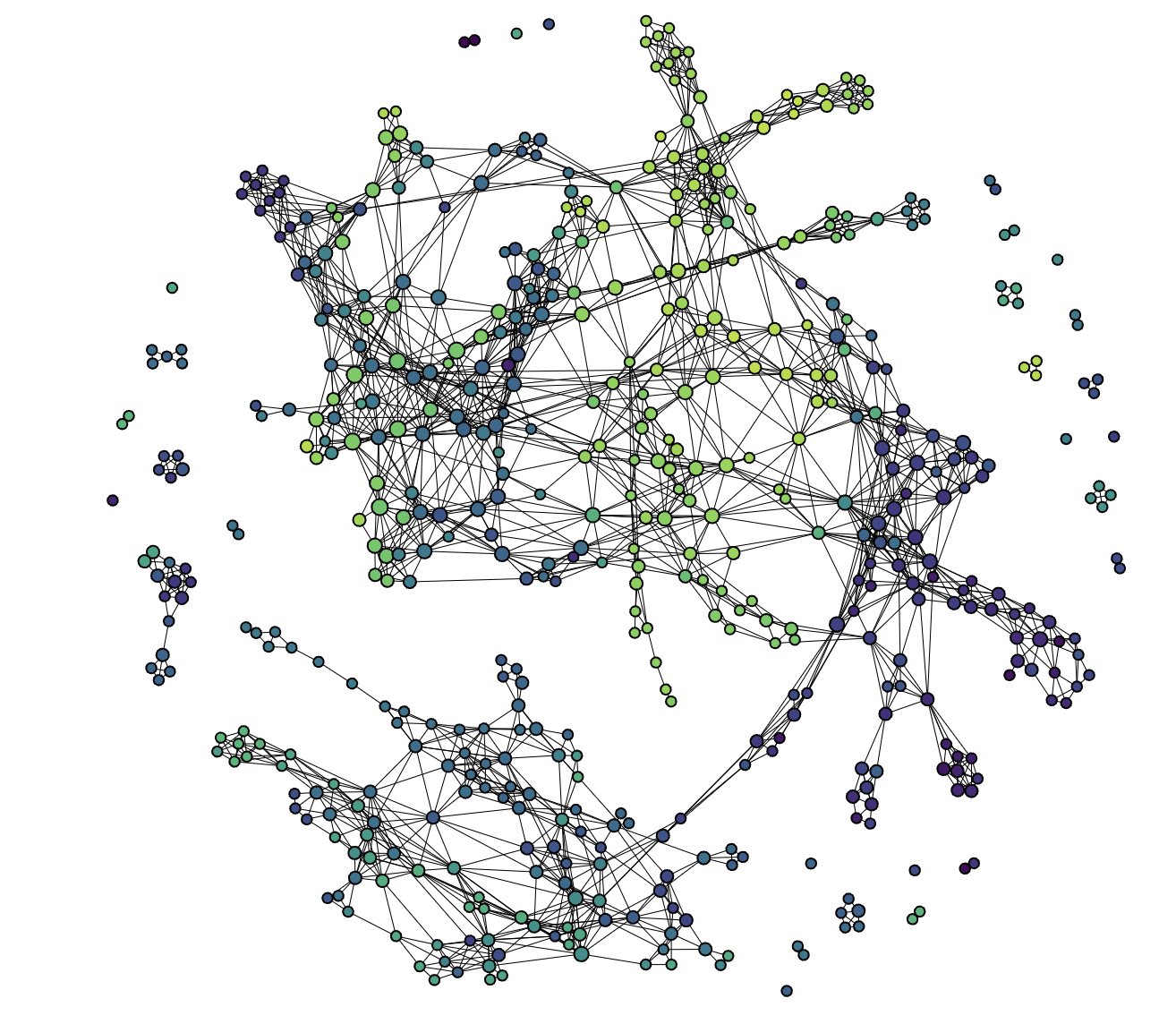}\vspace{-11pt}
                \caption{Midwest 3D (with the dates)}\label{fig:midwest_3D}
	\end{subfigure}
        \vskip3mm
	\caption{Midwest Clustering via Mapper (with the dates)}\label{fig:midwest_data}
\end{figure}

\begin{figure}[!htbp]
	\begin{subfigure}[b]{0.4\textwidth}
		\includegraphics[width=\textwidth]{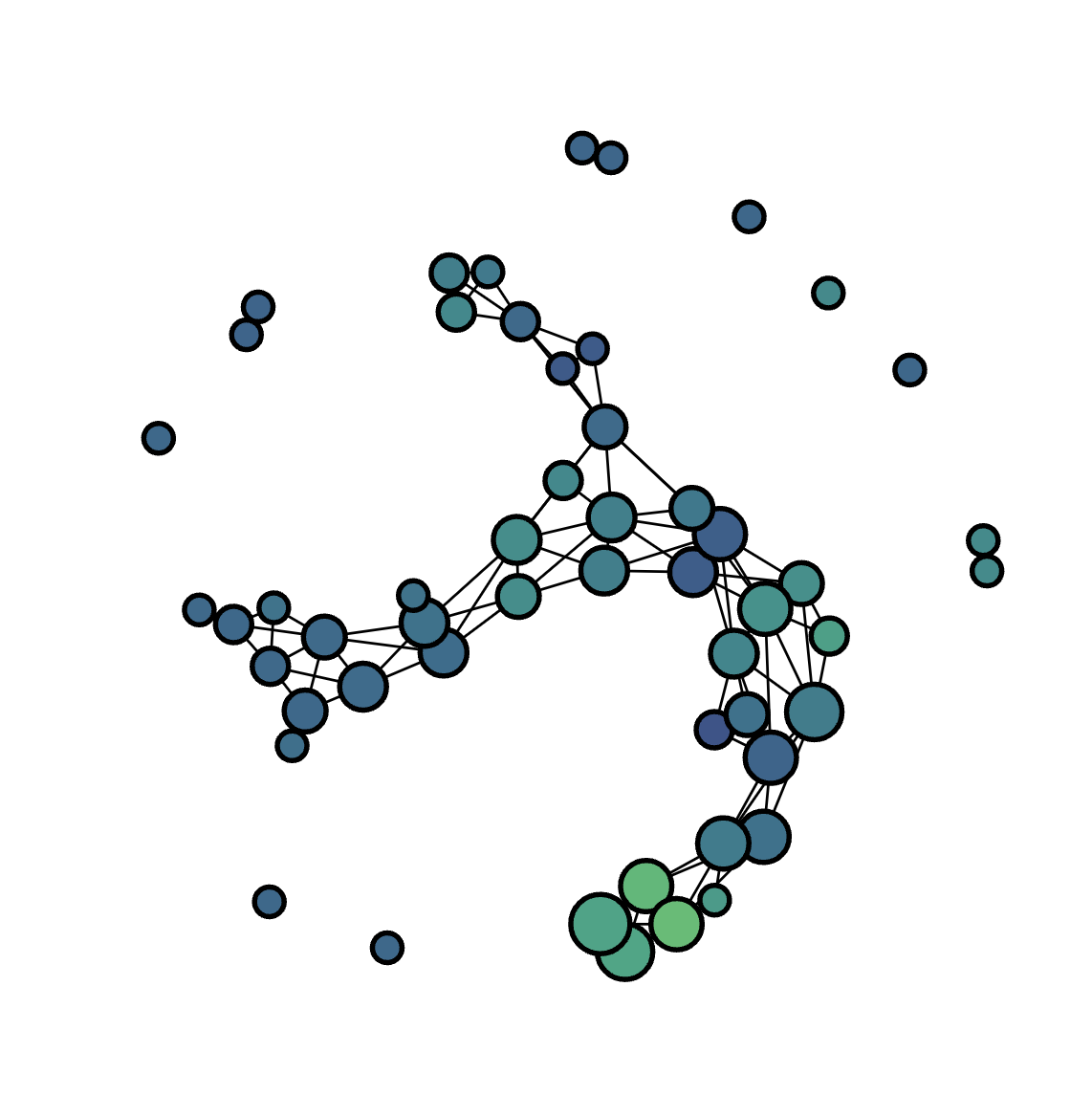}\vspace{-11pt}
                \caption{Northeast 2D (without the dates)}\label{fig:northeast_2D_noDates}
	\end{subfigure}
	\begin{subfigure}[b]{0.4\textwidth}
		 \includegraphics[width=\textwidth]{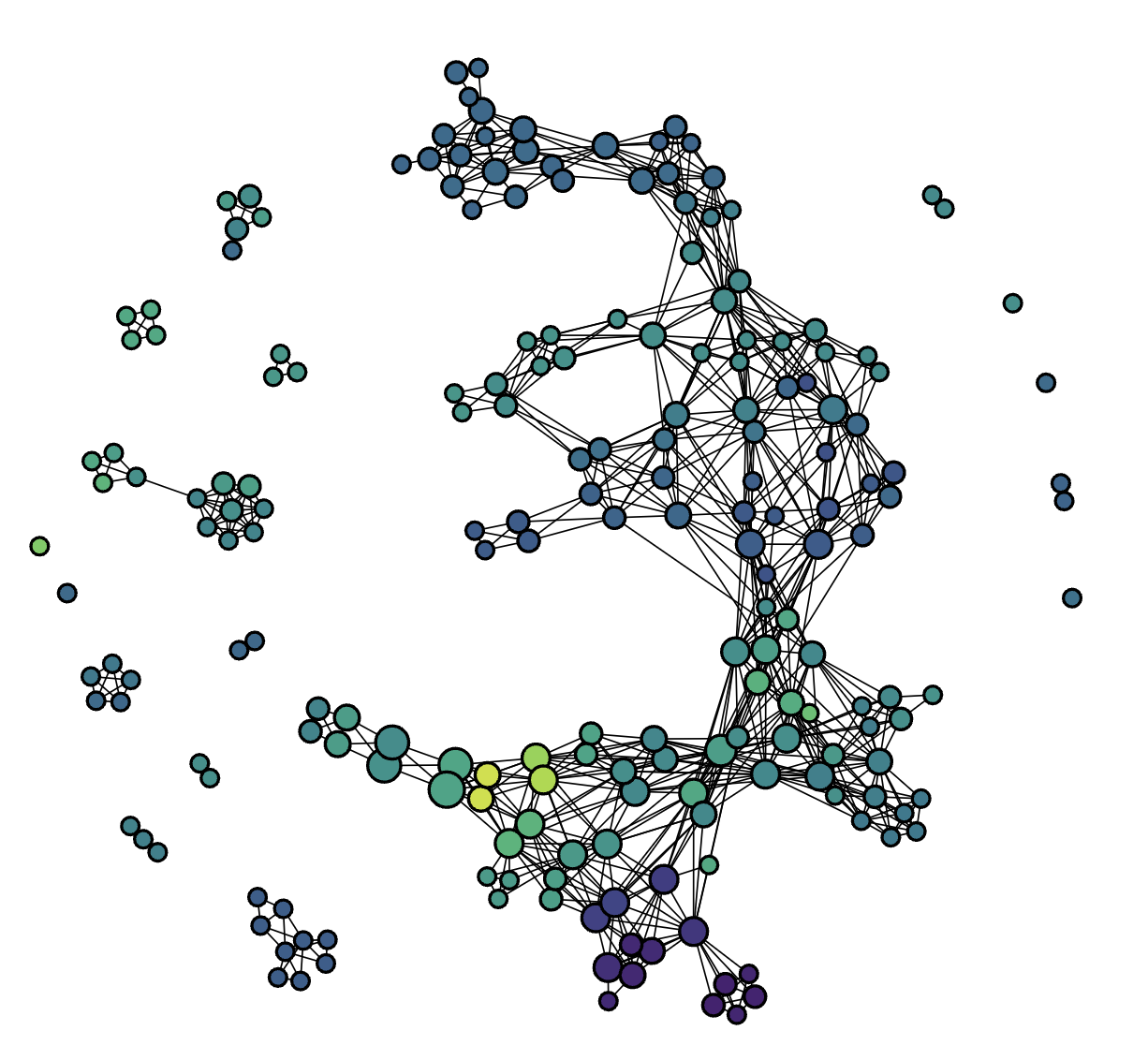}\vspace{-11pt}
                \caption{Northeast 3D (without the dates)}\label{fig:northeast_3D_noDates}
	\end{subfigure}
    
	\vskip3mm

        \begin{subfigure}[b]{0.38\textwidth}
		\includegraphics[width=\textwidth]{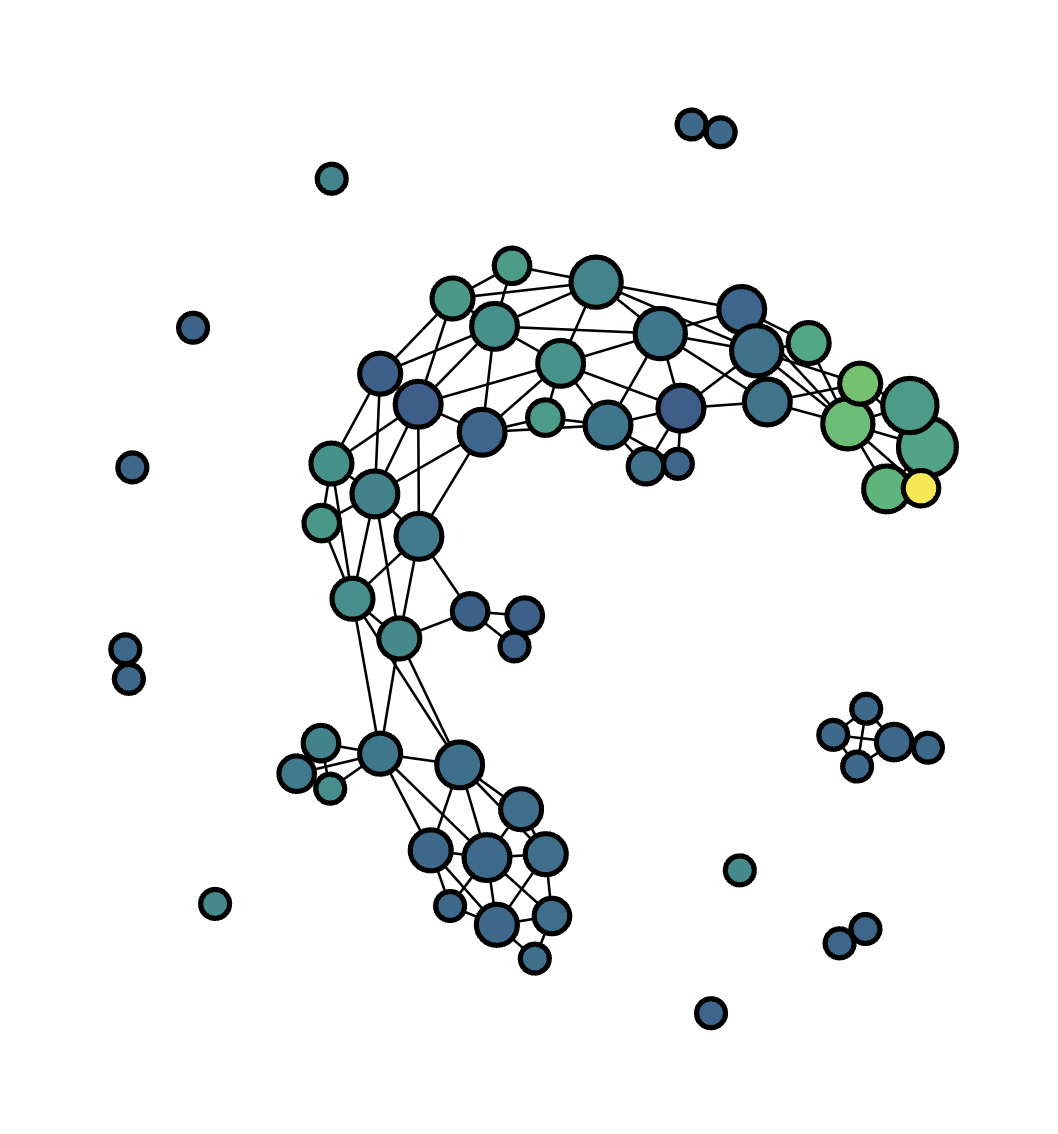}\vspace{-11pt}
                \caption{Northeast 2D (with the dates)}\label{fig:noreast_2D}
	\end{subfigure}
	~
	\begin{subfigure}[b]{0.42\textwidth}
		\includegraphics[width=\textwidth]{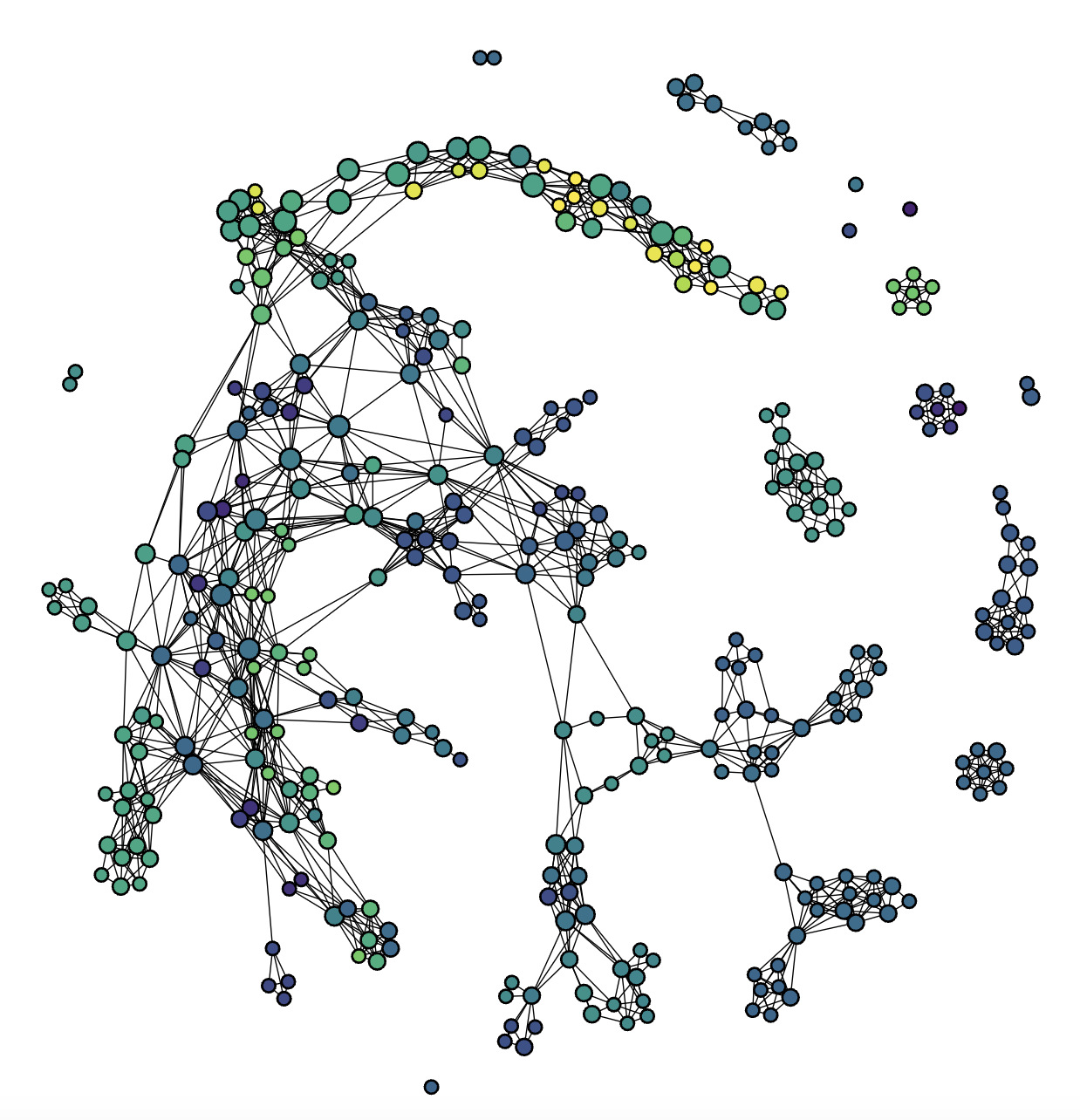}\vspace{-11pt}
                \caption{Northeast 3D (with the dates)}\label{fig:northeast_3D}
	\end{subfigure}
        \vskip3mm
	\caption{Northeast Clustering via Mapper}\label{fig:northeast_data}
\end{figure}

\begin{figure}[!htbp]
	\begin{subfigure}[b]{0.41\textwidth}
		\includegraphics[width=\textwidth]{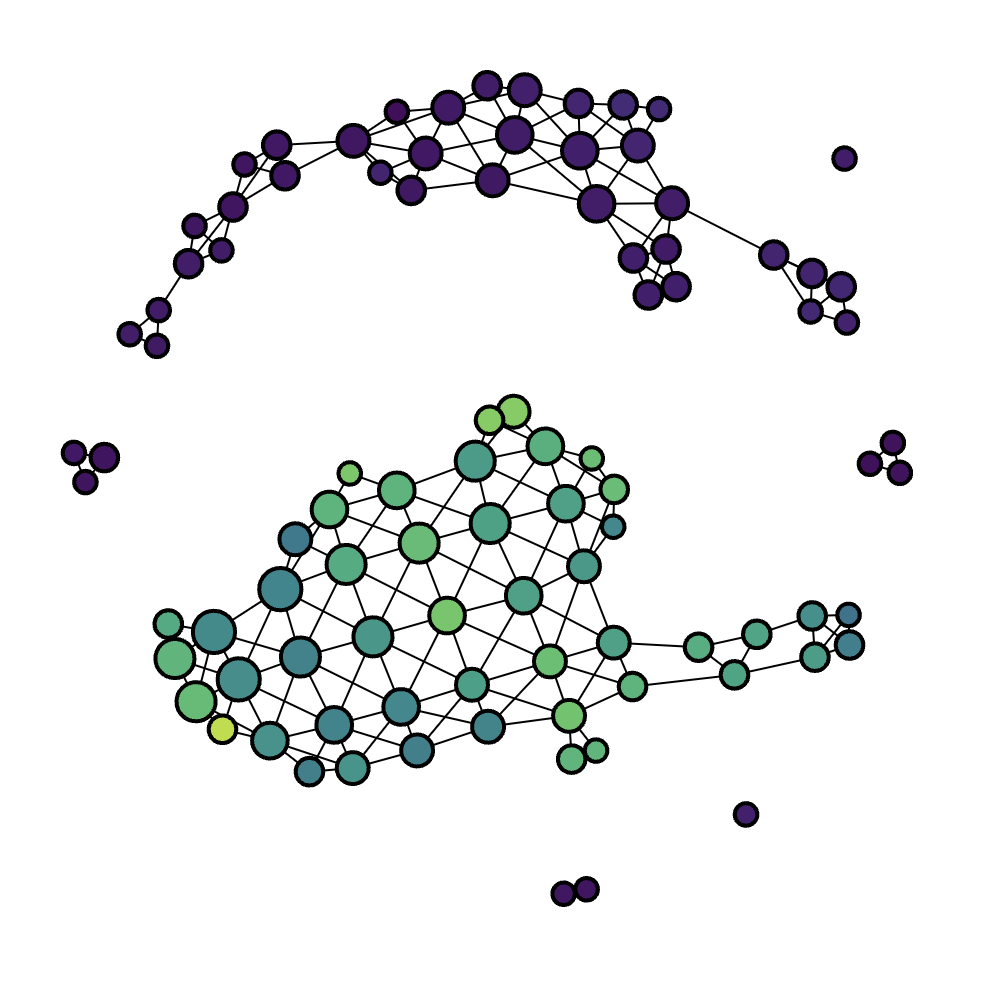}\vspace{-11pt}
                \caption{South 2D (without the dates)} \label{fig:south_2D_noDates}
	\end{subfigure}
	\begin{subfigure}[b]{0.41\textwidth}
		 \includegraphics[width=\textwidth]{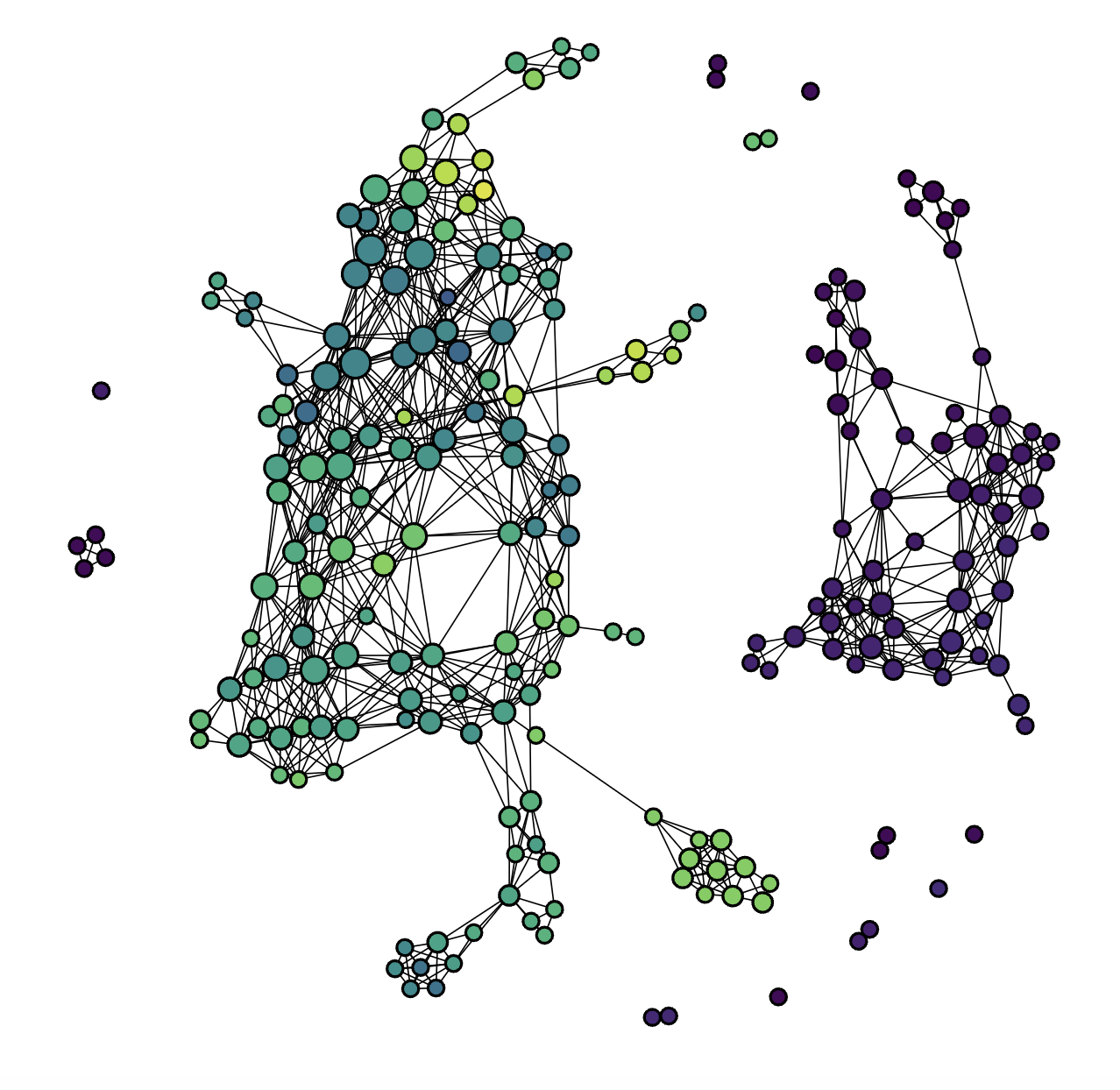}\vspace{-11pt}
                \label{fig:south_3D_noDates}
                \caption{South 3D (without the dates)}
	\end{subfigure}
    
	\vskip3mm

    \begin{subfigure}[b]{0.35\textwidth}
	    \includegraphics[width=\textwidth]{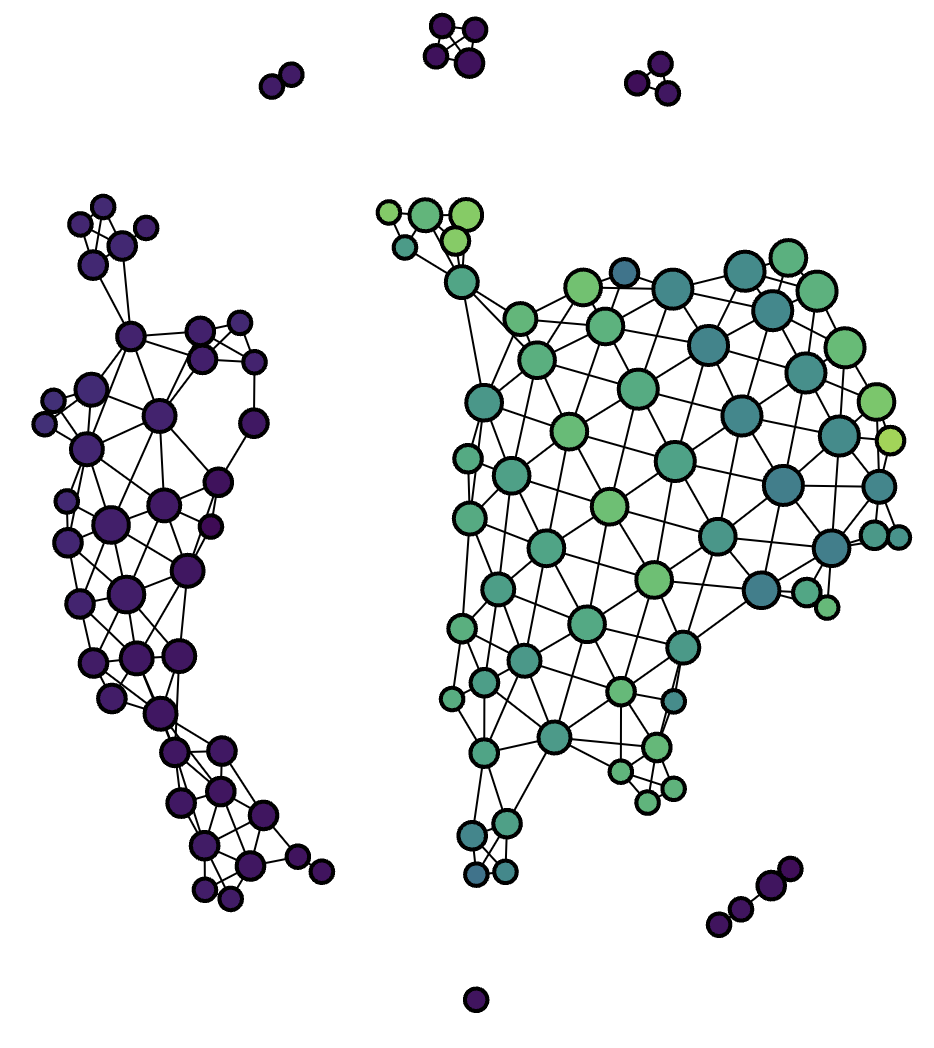}\vspace{-11pt}
                \label{fig:south_2D}
                \caption{South 2D (with the dates)}
	\end{subfigure}
	~
	\begin{subfigure}[b]{0.48\textwidth}
		\includegraphics[width=\textwidth]{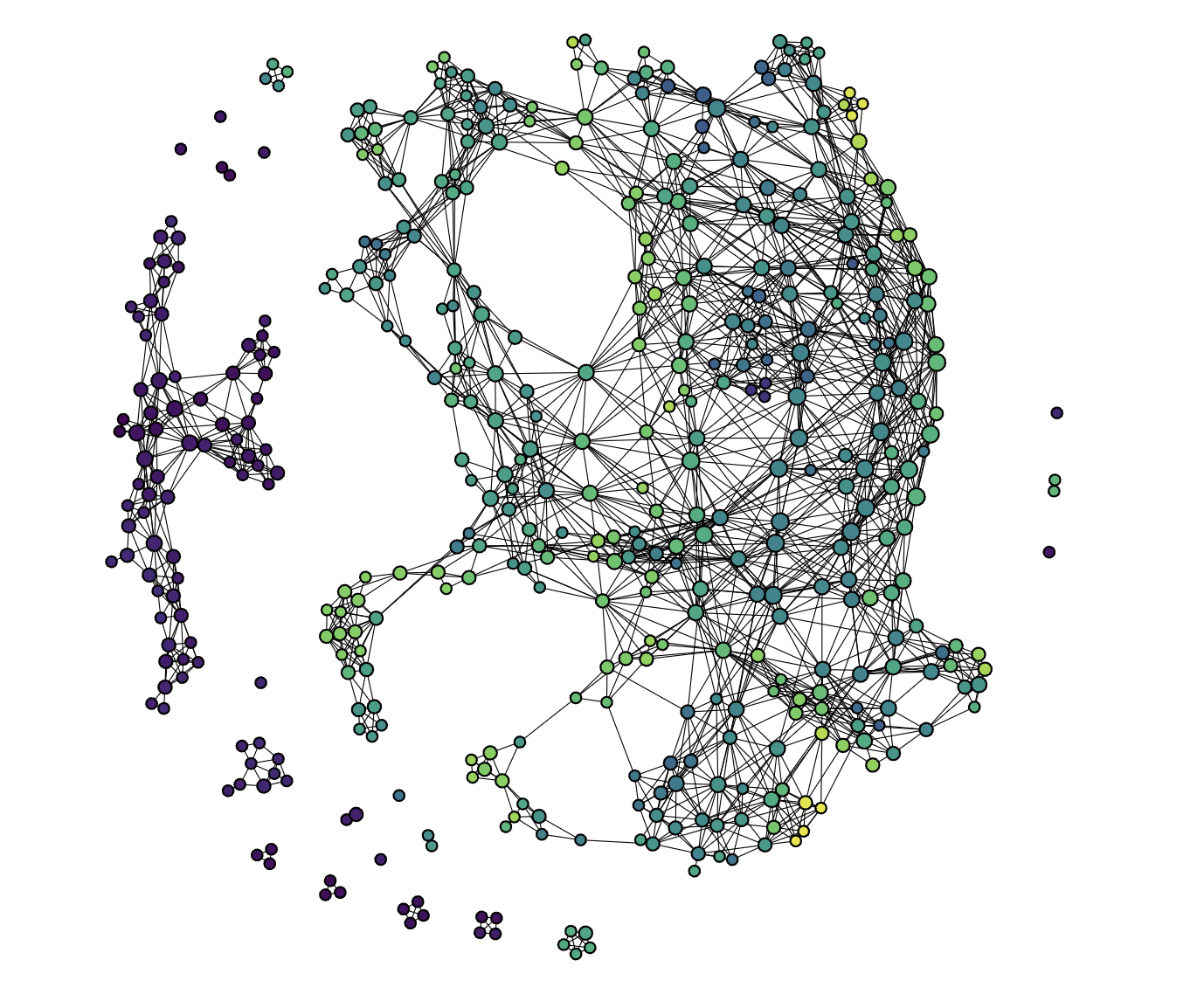}\vspace{-11pt}
                \label{fig:south_3D}
                \caption{South 3D (with the dates)}
	\end{subfigure}
        \vskip3mm
	\caption{South Clustering via Mapper}\label{fig:south_data}
\end{figure}

\begin{figure}[!htbp]
	\begin{subfigure}[b]{0.4\textwidth}
		\includegraphics[width=\textwidth]{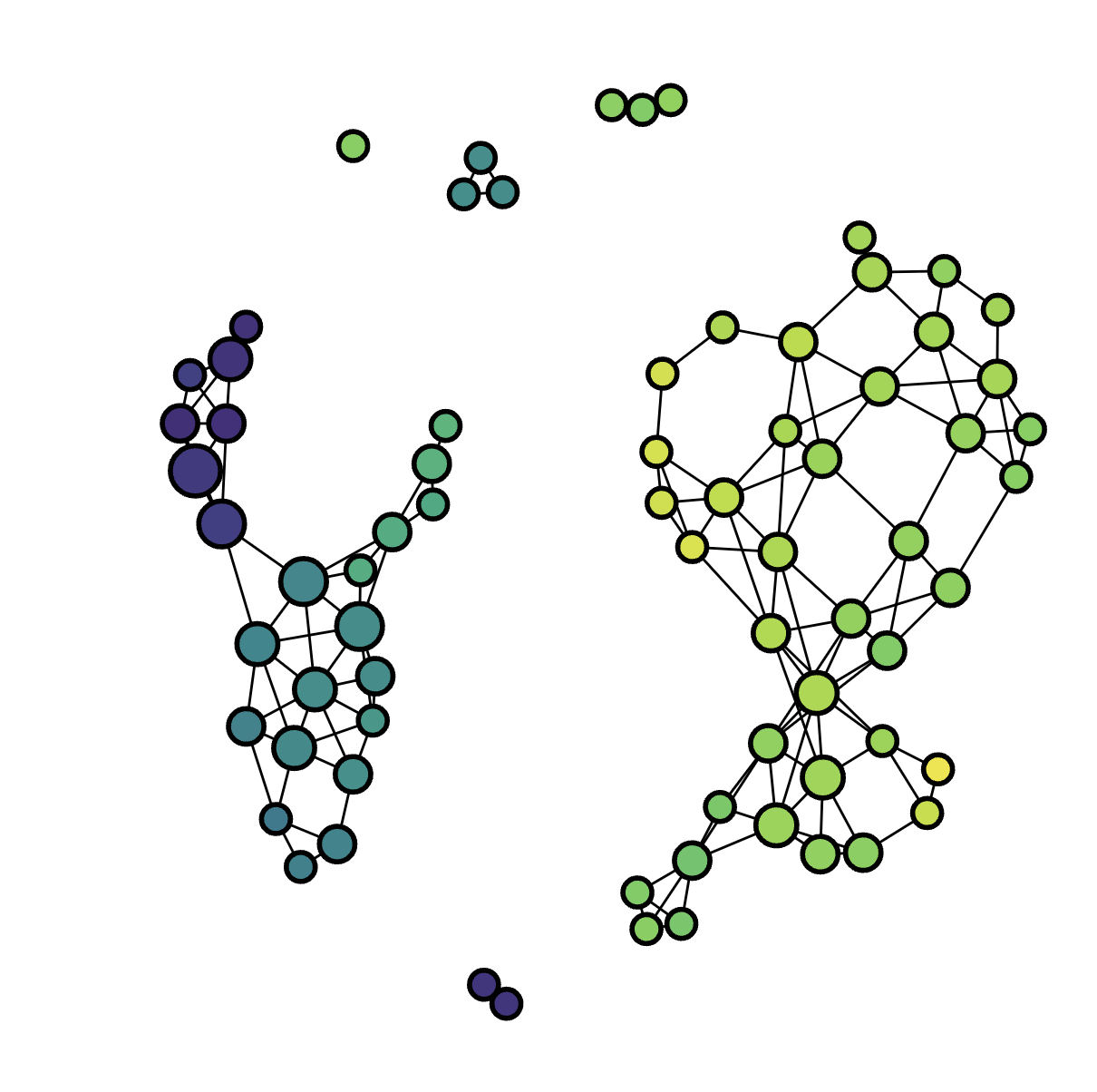}\vspace{-11pt}
                \caption{Non-Contiguous US 2D (without the dates)} \label{fig:outliers_2D_noDates}
	\end{subfigure}
	\begin{subfigure}[b]{0.48\textwidth}
		 \includegraphics[width=\textwidth]{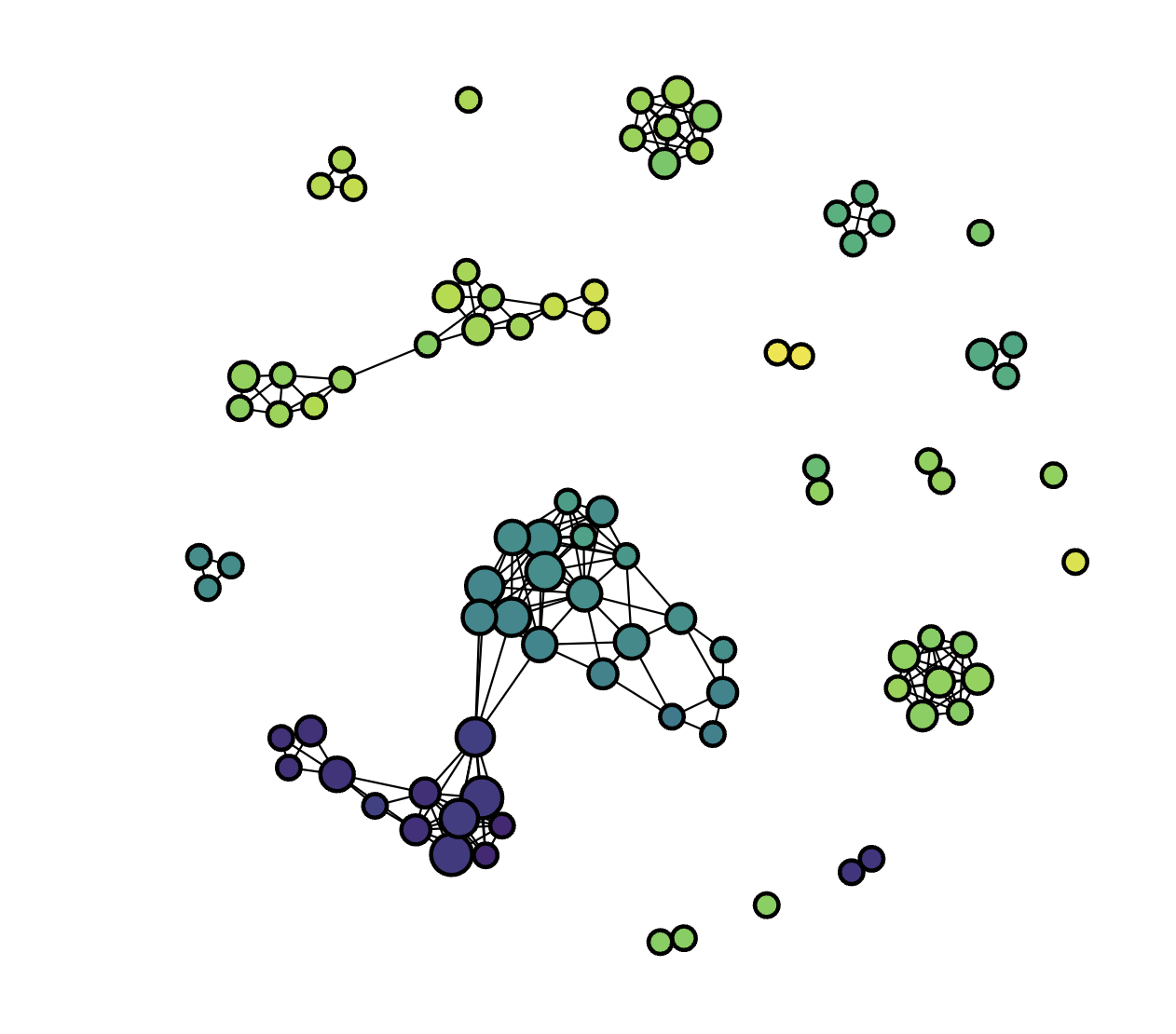}\vspace{-11pt}
                \caption{Non-Contiguous US 3D (without the dates)}\label{fig:outliers_3D_noDates}
	\end{subfigure}
    
	\vskip3mm

        \begin{subfigure}[b]{0.45\textwidth}
		\includegraphics[width=\textwidth]{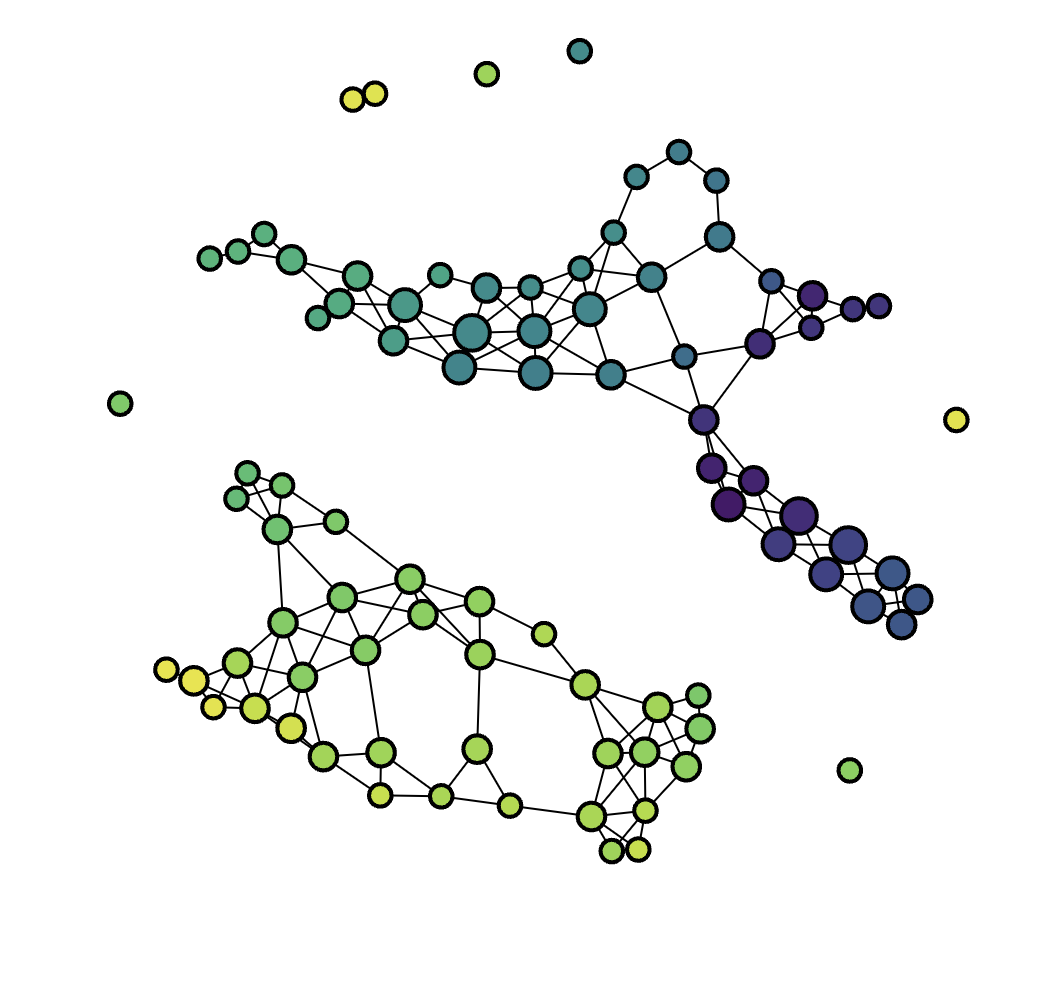}\vspace{-11pt}
                \caption{Non-Contiguous US 2D (with the dates)}\label{fig:outliers_2D}
	\end{subfigure}
	~
	\begin{subfigure}[b]{0.48\textwidth}
		\includegraphics[width=\textwidth]{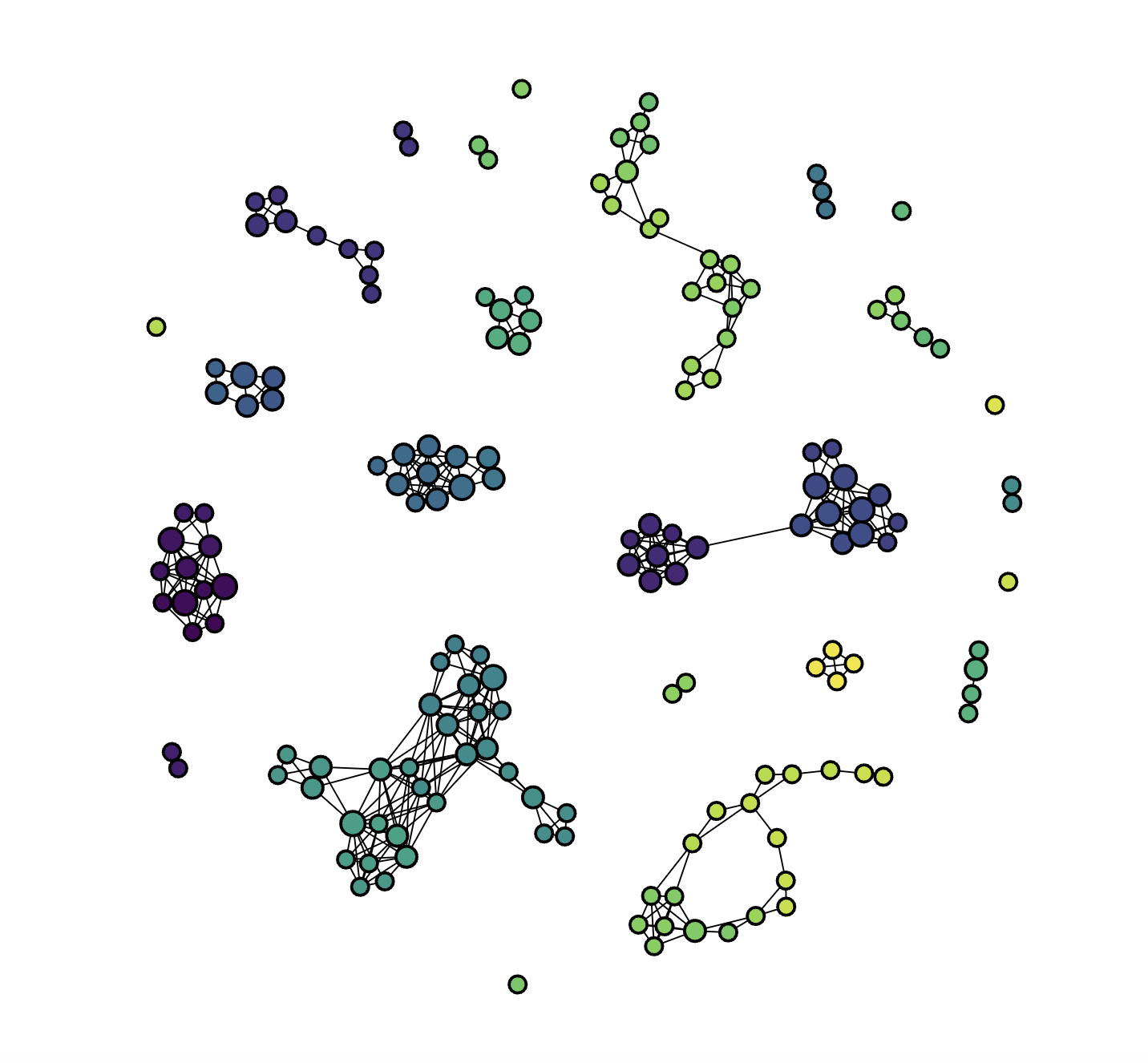}\vspace{-11pt}
                \caption{Non-Contiguous US 3D (with the dates)}\label{fig:outliers_3D}
	\end{subfigure}
        \vskip3mm
	\caption{Non-Contiguous US Clustering via Mapper}\label{fig:outliers_data}
\end{figure}
\end{document}